\documentclass[11pt,twoside,a4paper]{article}
\usepackage{amsmath,amssymb,mathrsfs}
\usepackage{amsthm}
\usepackage[utf8x]{inputenc}
\usepackage{bm}
\usepackage{avant}
\usepackage[english]{babel}
\usepackage{thmtools,thm-restate}
\usepackage[nottoc]{tocbibind}
\usepackage{hyperref}
\usepackage{graphicx}
\usepackage{enumitem}
\usepackage{tikz}
\usepackage{bbm}
\usepackage[inner=2.4cm,outer=2.4cm,top=2.8cm,bottom=2.8cm]{geometry}
\usepackage{mathtools}
\usepackage{makeidx}
\usepackage{bold-extra}

\DeclareMathOperator{\id}{id}

\newcommand{\modu}[1]{\left |#1\right |}
\newcommand{\Leb}{\mathscr{L}}
\newcommand{\setR}{\mathbb{R}}
\newcommand{\setN}{\mathbb{N}}
\newcommand{\R}{\mathbb{R}}

\newcommand{\p}{\mathtt p} %projection
\newcommand{\de}{\ensuremath{\, \mathrm d}} % in the integrals

\newcommand{\suchthat}{\ensuremath{\,:\,}} % such that inside, for example, the sets definition
\newcommand\restr[2]{{% we make the whole thing an ordinary symbol
  \left.\kern-\nulldelimiterspace % automatically resize the bar with \right
  #1 % the function
  \right|_{#2} % this is the delimiter
  }}

\newcommand{\weakto}{\rightharpoonup}

\newcommand{\CD}{\mathsf{CD}}

\newcommand{\MCP}{\mathsf{MCP}}

\newcommand{\di}{\mathsf d} %standard mms distance notation
\newcommand{\m}{\mathfrak m} %standard mms measure notation 
\DeclareMathOperator{\Ent}{Ent}
\DeclareMathOperator{\Geo}{Geo}
\DeclareMathOperator{\OptGeo}{OptGeo}
\DeclareMathOperator{\OptPlans}{OptPlans}
\newcommand{\Prob}{\mathscr{P}}
\newcommand{\ProbTwo}{\mathscr{P}_2}
\newcommand{\ppi}{{\mbox{\boldmath$\pi$}}}

\usepackage{titlesec}
\titleformat{\section}
  {\scshape\large\centering}
  {\thesection}{0.5em}{}
\titleformat{\subsection}
  {\scshape\centering}
  {\thesubsection}{0.4em}{}  

\author{Mattia Magnabosco}
\title{\textbf{Example of an Highly Branching CD Space}}
\date{}

\newtheoremstyle{remark}% <name>
        {10pt}% <Space above>
        {10pt}% <Space below>
        {}% <Body font>
        {}% <Indent amount>
        {\itshape}% <Theorem head font>
        {.}% <Punctuation after theorem head>
        {.4em}% <Space after theorem head>
        {}% <Theorem head spec (can be left empty, meaning 'normal')>

\newtheoremstyle{proof}% <name>
        {10pt}% <Space above>
        {10pt}% <Space below>
        {}% <Body font>
        {}% <Indent amount>
        {\itshape}% <Theorem head font>
        {.}% <Punctuation after theorem head>
        {.4em}% <Space after theorem head>
        {}% <Theorem head spec (can be left empty, meaning 'normal')>
        
\newtheoremstyle{definition}% <name>
        {10pt}% <Space above>
        {10pt}% <Space below>
        {}% <Body font>
        {}% <Indent amount>
        {\bfseries}% <Theorem head font>
        {.}% <Punctuation after theorem head>
        {.4em}% <Space after theorem head>
        {}% <Theorem head spec (can be left empty, meaning 'normal')>        

\newtheoremstyle{theorem}% <name>
        {10pt}% <Space above>
        {10pt}% <Space below>
        {\slshape}% <Body font>
        {}% <Indent amount>
        {\bfseries}% <Theorem head font>
        {.}% <Punctuation after theorem head>
        {.4em}% <Space after theorem head>
        {}% <Theorem head spec (can be left empty, meaning 'normal')>

\theoremstyle{theorem}
\newtheorem{theorem}{Theorem}[section]
\newtheorem{prop}[theorem]{Proposition}
\newtheorem{corollary}[theorem]{Corollary}
\newtheorem{lemma}[theorem]{Lemma}

\theoremstyle{definition}
\newtheorem{definition}[theorem]{Definition}

\theoremstyle{remark}

\theoremstyle{proof}
\newtheorem*{pro}{Proof}
\newenvironment{pr}{\begin{pro}%
 \pushQED{\qed}}%
 {\popQED\end{pro}}

\makeindex
\begin{document}
\maketitle

\begin{abstract}
   In \cite{ketterer2014failure} Ketterer and Rajala showed an example of metric measure space, satisfying the measure contraction property $\MCP(0,3)$, that has different topological dimensions at different regions of the space. In this article I propose a refinement of that example, which satisfies the $\CD(0,\infty)$ condition, proving the non-constancy of topological dimension for CD spaces. This example also shows that the weak curvature dimension bound, in the sense of Lott-Sturm-Villani, is not sufficient to deduce any reasonable non-branching condition. Moreover, it allows to answer to some open question proposed by Schultz in \cite{schultz2017existence}, about strict curvature dimension bounds and their stability with respect to the measured Gromov Hausdorff convergence. 
\end{abstract}

In their remarkable works Lott, Villani \cite{lottvillani} and Sturm \cite{sturm2006,sturm2006ii} introduced a weak notion of curvature dimension bounds, which strongly relies on the theory of Optimal Transport. Inspired by some results that hold in the Riemannian case, they defined a consistent notion of curvature dimension bound for metric measure spaces, that is known as CD condition. The metric measure spaces satisfying the CD condition are called CD spaces and enjoy some remarkable analytic and geometric properties.

In this work I present an example of an highly branching CD space, that shows how the weak curvature dimension bound is not sufficient to deduce any type of non-branching condition. This example is a refinement of the one by Ketterer and Rajala \cite{ketterer2014failure}, in particular, the topological and metric structure of the space will be essentially the same, while the reference measure will be more complicated. The metric measure space considered by Ketterer and Rajala in \cite{ketterer2014failure} satisfies the so called \textit{measure contraction property} $\MCP(0,3)$, but, as they pointed out in their work, it does not satisfy any CD condition. As I will highlight in the last section, the modification I am going to present allows to extend the some of the results in \cite{ketterer2014failure} to CD spaces. Moreover it gives an answer to some questions proposed by Schultz in \cite{schultz2017existence}, regarding the very strict CD condition and the stability of strict CD condition. 

Before going on, I want to emphasize that the CD condition is much more complicated to prove than the MCP condition. In fact there is much more freedom in choosing the marginals, and consequently it is much more difficult to avoid mass overlap in a Wasserstein geodesic. For this reason, not only will the proof be longer than Ketterer and Rajala's one, but it will also be more complicated. In particular, I will prove the CD condition for the desired metric measure space, by proving it for a sequence of metric measure spaces that converges to it, and using the stability result by Villani (Theorem \ref{thm:stabilitycompact}). I will also take great inspiration by the work of Rajala in \cite{rajala2013failure}, which shows an efficient way to prove the CD condition in branching metric measure spaces, through the application of Jacobi equation.   
\\

\noindent The main consequences of the example I will present in this article are the following:
\begin{itemize}
    \item any reasonable non-branching condition does not hold in general for (weak) CD spaces,
    \item the existence of an optimal transport map between two absolutely continuous marginals is not granted in (weak) CD spaces, without assuming a non-branching condition,
    \item the very strict CD condition studied by Schultz (\cite{schultz2017existence},\cite{Schultz2019EquivalentDO},\cite{Schultz2019OnOO}) is strictly stronger than the weak one,
    \item the constancy of topological dimension does not hold in general for (weak) CD spaces,
    \item the strict CD condition I will define in the last section is not stable with respect of the measured Gromov Hausdorff convergence,
    \item the strict CD condition is strictly stronger than the weak one.
\end{itemize}

Let me now briefly explain the structure of this work. In the first section I recall some preliminary results regarding both the basis of the Optimal Transport theory and CD spaces. In the second section I am going to introduce the metric measure spaces which will be the subject of the rest of the article. In section \ref{section:algebraiclemmas} I simply state and prove some algebraic lemma, which will help me doing the subsequent computations. In the fourth section I present some result by Rajala \cite{rajala2013failure}, related to Jacobi equation and how it can be used to prove entropy convexity. This theory requires the existence of a suitable midpoint selection map, which will be introduced in Section \ref{sec:defimidpoint}, where I will also prove some of its properties. The sixth section contains the proof of the main theorem, that puts together all the results proven so far. The last section aims to draw all the conclusions listed before.

\section{Preliminary Notions}\label{chapter:introduction}

This first section aims to collect all the preliminary results this work needs in order to be self contained. In particular I am going to introduce the Wasserstein space and the entropy functional on it, being then able to define the notions of curvature dimension bound and CD space. Moreover I am going to briefly discuss the relation between curvature dimension bound and non-branching conditions, that is one of the main motivation for this work. Finally I will define the measured Gromov Hausdorff convergence of metric measure spaces, stating in the end the stability of curvature dimension bounds with respect to this convergence.

\subsection{The Wasserstein Space}

%EDIT intro rapida
Denote by $\Prob(X)$ the set of Borel probability measures on a Polish metric space $(X,\di)$. Given two measures $\mu,\nu\in\Prob(X)$ and a Borel cost function $c:X\times X \to [0,\infty]$, the Optimal Transport problem asks to find minima and minimizers of the quantity
\begin{equation}\label{eq:OTproblem}
\min  \int_{X\times Y} c(x,y) \de \pi(x,y)  ,
\end{equation}
where $\pi$ varies among all probability measures in $\Prob(X\times X)$ with first marginal equal to $\mu$ and second marginal equal to $\nu$. If the cost function $c$ is lower semicontinuous, 
the minimum in \eqref{eq:OTproblem} is attained.
The minimizers of this problem are called optimal transport plans and the set of all of them will be denoted by $\OptPlans(\mu,\nu)$. An optimal transport plan $\pi\in\OptPlans(\mu,\nu) $ is said to be induced by a map if there exists a $\mu$-measurable map $T:X \to X$ so that $\pi=(\id,T)_\# \mu$, such a map $T$ will be called optimal transport map.

\noindent A fundamental approach in facing the Optimal Transport problem is the one of $c$-duality, which allows to prove some very interesting and useful results. Below I report only the most basic statement, which is the only result I will need in this work.

\begin{definition}
A set $\Gamma\subset X\times X$ is said to be $c$-cyclically monotone if 
\begin{equation*}
    \sum_{i=1}^{N} c\left(x_{i}, y_{\sigma(i)}\right) \geq \sum_{i=1}^{N} c\left(x_{i}, y_{i}\right)
\end{equation*}
for every $N\geq1$, every permutation $\sigma$ of $\{1,\dots,N\}$ and every $(x_i,y_i)\in \Gamma$ for $i=1,\dots,N$.
\end{definition}

\begin{prop}\label{prop:cmonotonicity}
Let $X$ be a Polish space and $c:X\times X \to [0,\infty]$ a lower semicontinuous cost function. Then every optimal transport plan $\pi\in \OptPlans(\mu,\nu)$ such that $\int c \de \pi<\infty$ is concentrated in a $c$-cyclically monotone set. 
\end{prop}

From now on I am going to consider the Optimal Transport problem in the special case in which the cost function is equal to the distance squared, that is $c(x,y)=\di^2(x,y)$. In this context the minimization problem induces the so called Wasserstein distance on the space $\ProbTwo(X)$ of probabilities with finite second order moment, that is

\begin{equation*}
    \ProbTwo(X):= \left\{\mu\in \Prob(X) \suchthat \int \di^2(x,x_0) \de\mu(x) <\infty \text{ for one (and thus all) } x_0\in X\right\}.
\end{equation*}

\begin{definition}[Wasserstein distance]
Given two measures $\mu,\nu\in \ProbTwo(X)$ define their Wasserstein distance $W_2(\mu,\nu)$ as 
\begin{equation*}
W_2^2(\mu,\nu) := \min \left\{ \int d^2(x,y) \de \pi(x,y) \suchthat \pi \in \Gamma(\mu,\nu) \right \}.
\end{equation*}
\end{definition}

\noindent It is easy to realize that $W_2$ is actually a distance on $\ProbTwo(X)$, moreover $(\ProbTwo(X),W_2)$ is a Polish metric space.

\noindent Let me now deal with the geodesic structure of $(\ProbTwo(X),W_2)$, which, as the following statement shows, is heavily related to the one of the base space $(X,\di)$. First of all, notice that every measure $\ppi\in\Prob(C([0,1],X))$ induces a curve $[0,1]\ni t \to \mu_t=(e_t)_\# \ppi \in \Prob(X)$, therefore in the following I will consider measures in $\Prob(C([0,1],X))$ in order to consider curves in the Wasserstein space.

\begin{prop}\label{prop:optgeo}
If $(X,\di)$ is a geodesic space than $(\ProbTwo(X),W_2)$ is geodesic as well.
In particular, given two measures $\mu,\nu\in \ProbTwo(X)$, the measure $\ppi\in\Prob(C([0,1],X))$ is a constant speed Wassertein geodesic connecting $\mu$ and $\nu$ if and only if it is concentrated in $\Geo(X)$ (that is the space of constant speed geodesics in $(X,\di)$) and $(e_0,e_1)_\#\ppi\in\OptPlans(\mu,\nu)$. In this case it is said that $\ppi$ is an optimal geodesic plan between $\mu$ and $\nu$ and this will be denoted as $\ppi\in\OptGeo(\mu,\nu)$.
\end{prop}

\subsection{Curvature Dimension Bounds}

In this subsection I introduce the notions of curvature dimension bound and CD space, pioneered by Lott and Villani \cite{lottvillani} and Sturm \cite{sturm2006,sturm2006ii}. Their definition relies on the notion of entropy functional. As it will be soon clear, the most appropriate framework in which deal with the entropy functional, is the one of metric measure spaces.

\begin{definition}
A metric measure space is a triple $(X,\mathsf{d},\mathfrak{m})$, where $(X,\mathsf{d})$ is a Polish metric space and $\mathfrak{m}$ is a non-negative and non-null Borel measure on $X$, finite on bounded sets.
\end{definition}

\noindent In this work I will only deal with compact metric measure spaces, and in particular $\m(X)<\infty$. Let me now properly define the entropy functional.

\begin{definition}
In a metric measure space $(X,\di,\m)$, define the relative entropy functional with respect to the reference measure $\m$ $\Ent:\ProbTwo(X)\to \R \cup \{+\infty\}$ as 
\begin{equation*}
    \Ent(\mu):= 
    \begin{cases}
    \int \rho \log \rho \de \m  &\text{if }\mu\ll\nu \text{ and } \mu=\rho\m\\
    +\infty &\text{otherwise}
    \end{cases}.
\end{equation*}
\end{definition}

\noindent In this context of this work the entropy functional $\Ent$ is lower semicontinuous with respect to the Wasserstein convergence, this is not always true in the non-compact case, when it might happen that $\m(X)=+\infty$. \\
I can now give the definitions of CD condition and CD space.

\begin{definition}\label{def:cdinfty}
A metric measure space $(X,\di,\m)$ is said to satisfy the (weak) $\CD(K,\infty)$ condition and to be a (weak) $\CD(K,\infty)$ space, if for every absolutely continuous measures  $\mu_0,\mu_1\in\ProbTwo(X)$ there exists a Wasserstein geodesic with constant speed $(\mu_t)_{t\in[0,1]}\subset\ProbTwo(X)$ connecting them, along which the relative entropy functional is $K$-convex, that is 
\begin{equation}\label{eq:CD(K,infty)}
    \Ent(\mu_t) \leq (1-t) \Ent(\mu_0) + t\Ent(\mu_1) -t(1-t)\frac K2  W_2^2(\mu_0,\mu_1), \qquad \text{for every }t\in [0,1].
\end{equation}
Moreover $(X,\mathsf{d},\mathfrak{m})$ is said to satisfy a strong $\CD(K,\infty)$ condition and to be a strong $\CD(K,\infty)$ space if, for every absolutely continuous measures $\mu_0,\mu_1\in\ProbTwo(X)$, the relative entropy functional is $K$-convex along every Wasserstein geodesic with constant speed connecting them. 
\end{definition}

\noindent Let me also state a very useful proposition, which provides a simple strategy to prove the (weak) $\CD(K,\infty)$ condition. Its proof can be found in \cite{sturm2006}, anyway I present a brief sketch of it, in order to be self contained.

\begin{prop}\label{prop:midpoint}
The metric measure space $(X,\di,\m)$ is a $\CD(K,\infty)$ space if for every pair of absolutely continuous measures $\mu_0,\mu_1\in\ProbTwo(X)$ there exists a midpoint $\eta\in\ProbTwo(X)$ of $\mu_0$ and $\mu_1$,
absolutely continuous with respect to $\m$, satisfying
\begin{equation}\label{eq:midpointconv}
    \Ent(\eta) \leq \frac12\Ent(\mu_0) + \frac12 \Ent(\mu_1) - \frac K8 W_2^2(\mu_0,\mu_1). 
\end{equation}
\end{prop}

\begin{pr}
Given two absolutely continuous measures $\mu_0,\mu_1\in\ProbTwo(X)$, define $\mu_{1/2}$ as a midpoint of $\mu_0$ and $\mu_1$ satisfying \eqref{eq:midpointconv}. Similarly define $\mu_{1/4}$ as a midpoint of $\mu_0$ and $\mu_{1/2}$ and $\mu_{3/4}$ as a $\mu_{1/2}$ and $\mu_1$, both satisfying \eqref{eq:midpointconv}. Proceeding in this way, it is possible to define $\mu_t$ for all dyadic times $t\in \big\{\frac{k}{2^h}:h\in \setN^+, k=1,\dots, 2^h-1\big\}$. An easy induction argument on $h$ shows that
\begin{equation*}
    \Ent(\mu_t) \leq (1-t) \Ent(\mu_0) + t\Ent(\mu_1) -t(1-t)\frac K2  W_2^2(\mu_0,\mu_1).
\end{equation*}
 for every dyadic time $t\in \{\frac{k}{2^h}:h\in \setN^+, k=1,\dots, 2^h-1\}$. Defining the geodesic $(\mu_t)_{t\in[0,1]}$ as the continuous extension, the lower semicontinuity of the entropy ensures that, for every $t\in[0,1]$, the measure $\mu_t$ satisfies the equation \eqref{eq:CD(K,infty)}. 
\end{pr}

In the last part of this subsection I want to present the relation between curvature dimension bounds and non-branching conditions. 
The most important result in this context was proven by Rajala and Sturm in \cite{rajalasturm}:

\begin{theorem}\label{thm:RajalaSturm}
Every strong $\CD(K,\infty)$ metric measure space $(X,\di,\m)$ is essentially non-branching, that is for every absolutely continuous measures $\mu_0,\mu_1\in\ProbTwo(X)$, every optimal geodesic plan connecting them is concentrated on a non-branching set of geodesics.
\end{theorem} 

\noindent The work of Rajala and Sturm was then generalized by Schultz \cite{schultz2017existence} to the context of very strict CD spaces.

\begin{definition}\label{def:verystrict}
A metric measure space $(X,\mathsf{d},\mathfrak{m})$ is called a very strict $\CD(K,\infty)$ space
if for every absolutely continuous measures $\mu_0,\mu_1\in\ProbTwo(X)$
there exists an optimal geodesic plan $\eta\in\OptGeo(\mu_0,\mu_1)$, so that the entropy functional $\Ent$ satisfies the K-convexity inequality along $(\operatorname{restr}_{t_0}^{t_1})_\# (f\eta)$
for every $t_0<t_1\in [0,1]$, and for all bounded Borel functions $f : \Geo(X) \to \setR^+$ with $\int f \de \eta=1$.
\end{definition}

\noindent As the reader can easily notice, these spaces are not in general essentially non-branching, but they satisfy a weaker condition that I will call \textit{weak essentially non-branching}.

\begin{definition}[Weak Essentially Non-Branching]
A metric measure space $(X,\mathsf{d},\mathfrak{m})$ is said to be weakly essentially non-branching if for every absolutely continuous measures $\mu_0,\mu_1\in\ProbTwo(X)$, there exists an optimal geodesic plan connecting them, that is concentrated on a non-branching set of geodesics.
\end{definition}

\begin{theorem}[Schultz \cite{schultz2017existence}]\label{thm:schultz}
Every very strict $\CD(K,\infty)$ space is weakly essentially non-branching.
\end{theorem}

\noindent Notice that the very strict CD condition is intermediate between the weak and the strong one. It is easy to find examples of very strict CD spaces, which are not strong CD spaces, while it is not obvious if the very strict condition is strictly stronger than the weak one. In this work I am going to present an example of an highly branching (weak) CD space which is not very strict CD.

\subsection{Measured Gromov Hausdorff Convergence and Stability of CD Spaces}

In this subsection I introduce (following \cite{villani2008}) a notion of convergence for metric measure spaces, that is called measured Gromov Hausdorff convergence. Roughly speaking, it is the combination of Hausdorff topology for the metric side, and weak topology for the measure side. In order to properly define the measured Gromov Hausdorff convergence I have to preliminary introduce the notion of $\varepsilon$-isometry.

\begin{definition}
A measurable map $f:(X,\di,\m)\to(X',\di',\m')$ between two metric measure spaces is called an $\varepsilon$-isometry if 
\begin{enumerate}
    \item it almost preserves the distances, that is:
    \begin{equation*}
        \left|\di\left(f(x), f\left(x^{\prime}\right)\right)-\di\left(x, x^{\prime}\right)\right| \leq \varepsilon \quad \text{for every }x,x'\in X,
    \end{equation*}
    \item it is almost surjective, that is:
    \begin{equation*}
        \forall y \in X', \,\, \text{there exists } x \in X \text{ such that } \, \di(f(x), y) \leq \varepsilon.
    \end{equation*}
\end{enumerate}
\end{definition}

\begin{definition}
Let $(X_k,\di_k,\m_k)_{k\in\setN}$ and $(X,\di,\m)$ be compact metric measure spaces. It is said that the sequence $(X_k,\di_k,\m_k)_{k\in\setN}$ converges to $(X,\di,\m)$, in the measured Gromov Hausdorff sense, if for every $k$ there exists a measurable $\varepsilon_k$-isometry $f_k:X_k\to X$, where $\varepsilon_k\to0$, such that
\begin{equation}\label{eq:mGHconv}
    (f_k)_\# \m_k \weakto \m \quad \text{as }k\to \infty.
\end{equation}
\end{definition}

\noindent The measured Gromov Hausdorff convergence can be in some sense metrized by the $\mathbb{D}$ distance, introduced by Sturm in \cite{sturm2006ii}. Moreover in \cite{Gigli_2015} Gigli, Mondino and Savaré showed that some different notion of convergence for (pointed) metric measure spaces are equivalent to the (pointed) measured Gromov Hausdorff convergence.\\
I end this subsection, stating the stability of the (weak) CD condition with respect to the measured Gromov Hausdorff convergence.

\begin{theorem}\label{thm:stabilitycompact}
Let $(X_k,\di_k,\m_k)_{k\in\setN}$ be a sequence of compact metric measure spaces converging in the measured Gromov Hausdorff sense to a compact metric measure space $(X,\di,\m)$. Given $K\in\R$, if each $(X_k,\di_k,\m_k)_{k\in\setN}$ satisfies the weak curvature dimension condition $\CD(K,\infty)$, then also $(X,\di,\m)$ satisfies $\CD(K,\infty)$.
\end{theorem}

\section{Definition of the Metric Measure Spaces}

\begin{figure}
\begin{center}

\begin{tikzpicture} 

\filldraw[black!35!white] (-6,0)--(-1,0)--(-1,0.5)--(-6,0.5)--cycle;
\shade[left color=black!37!white,right color=black!5!white] (-1,0)--(6,0)--(6,3.5)--(1,1)..controls (0,0.5)..(-1,0.5)--cycle;
%\shade[top color=white,bottom color=black, opacity=0.35] (0,0)--(6,0)--(6,3)--(2,1)..controls (1,0.5)..(0,0.5)--cycle;
\draw[thick](-6,0)--(6,0);
\draw[thick](-6,0)--(-6,0.5);
\draw[thick](-6,0.5)--(-1,0.5);
\filldraw[black] (0,0) circle (1pt);
\draw[thick](6,0)-- (6,3.5)--(1,1);
\draw[thick](-1,0.5)..controls (0,0.5)..(1,1); 
\filldraw[white] (-7,0) circle (1pt);
\draw[<->,thick] (-5,0)--(-5,0.5);
\node at (-5.1,0.25)[label=east:\boldsymbol{$\varepsilon$}] {};
\node at (3,2.1)[label=north:$f_{k,\varepsilon}(x)$] {};
\end{tikzpicture}

\caption{The metric measure space $(X_{k,\varepsilon},\di_\infty, \mathfrak{m}_{k,K,\varepsilon})$ with $\varepsilon>0$.}
\label{fig:example}
\end{center}
\end{figure}
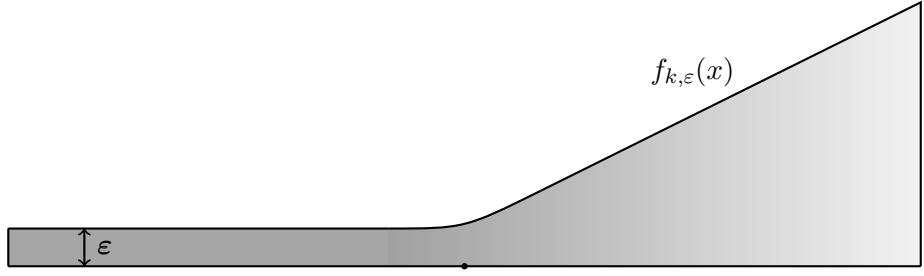

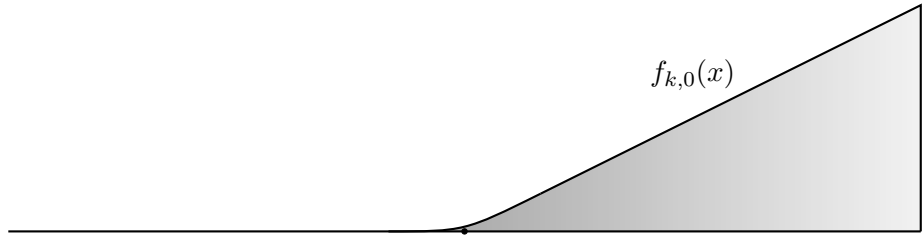
\begin{figure}[b]
\begin{center}

\begin{tikzpicture} 

\shade[left color=black!37!white,right color=black!5!white] (0,0)--(6,0)--(6,3)--(1,0.5)..controls (0,0)..(-1,0)--cycle;
\draw[thick](-6,0)--(6,0);
\filldraw[black] (0,0) circle (1pt);
\draw[thick](6,0)-- (6,3)--(1,0.5);
\draw[thick](-1,0)..controls (0,0)..(1,0.5); 
\filldraw[white] (-7,0) circle (1pt);
\node at (3,1.6)[label=north:$f_{k,0}(x)$] {};
\end{tikzpicture}

\caption{The metric measure space $(X_{k,0},\di_\infty, \mathfrak{m}_{k,K,0})$.}
\label{fig:example2}
\end{center}
\end{figure}

In this section I am going to introduce the metric measure spaces that will be studied in the rest of this work. The definitions that follow are actually involved and quite complicated, thus I invite the reader to look at Figure \ref{fig:example} and Figure \ref{fig:example2}, in order to better understand.

\begin{definition}
Fix $\varepsilon, k \in \setR$ such that $0\leq\varepsilon<k<\frac12$ and let $\varphi:\setR \to [0,1]$ be a continuous function such that $\int \varphi =k$, $\varphi\leq1$ on the whole $\R$ and $\varphi=0$ outside $[-k,k]$. Define the function $f_{k,\varepsilon}:[-1,1]\to \setR^+$ prescribing
\begin{equation*}
    f_{k,\varepsilon}''(x)= \varphi(x), \qquad f_{k,\varepsilon}'(-1)=0, \qquad f_{k,\varepsilon}(-1)=\varepsilon.
\end{equation*}
Consequently define the set $$X_{k,\varepsilon}=\{(x,y)\in \setR^2\suchthat x \in [-1,1] \text{ and } 0\leq y\leq f_{k,\varepsilon}(x)\}.$$
\end{definition}

\noindent In the following I will use this notation:
\begin{equation*}
    L= X_{k,0} \cap \big( \{f_{k,0}=0\} \times \R \big)\quad \text{and} \quad C= X_{k,0} \cap \big( \{f_{k,0}\ne0\} \times \R \big),
\end{equation*}
for sake of simplicity, I will not explicit the $k$ dependence, because it will be clear from the context. 

\begin{definition}
Given $\varepsilon, k \in \setR$ with $0<\varepsilon<k<\frac12$ and $K\geq1$, define the measure $\mathfrak{m}_{k,K,\varepsilon}$ on $X_{k,\varepsilon}$ as $$\mathfrak{m}_{k,K,\varepsilon}=m_{k,K,\varepsilon}(x,y)\cdot\Leb^2|_{X_{k,\varepsilon}}:=\frac{1}{f_{k,\varepsilon}(x)} \exp\left(-K \left(\frac{y}{f_{k,\varepsilon}(x)} \right)^2\right) \Leb^2|_{X_{k,\varepsilon}}.$$
While for every $0<k<\frac12$ and $K\geq1$, define the measure $\mathfrak{m}_{k,K,0}$ on $X_{k,0}$ as $$\mathfrak{m}_{k,K,0}:=\restr{  C_{K}\cdot\chi_{\{f_{k,0}(x)=0\}}\cdot\mathcal{H}^1}{\{y=0\}}  + \chi_{\{f_{k,0}(x)>0\}} \frac{1}{f_{k,0}(x)} \exp\left(-K \left(\frac{y}{f_{k,0}(x)} \right)^2\right)\cdot \Leb^2|_{X_{k,0}},$$ where
\begin{equation*}
    C_{K}= \int_0^1 e^{-K y^2} \de y.
\end{equation*}
\end{definition}

\noindent  Notice that a simple change of variable shows that 
\begin{equation*}
    (\p_x)_\# \mathfrak{m}_{k,K,0}= (\p_x)_\# \mathfrak{m}_{k,K,\varepsilon}=   C_{K}\cdot\chi_{\{-1\leq x\leq 1\}}\cdot\mathcal{H}^1,
\end{equation*}
for every suitable $k$, $K$ and $\varepsilon$. Moreover, since I have imposed $\varepsilon<k$, it is easy to realize that $f_{k,\varepsilon}(x)<3k$ for every $x \in [-1,1]$ (see Figure \ref{fig:example}).

In the following I am going to prove that for suitable $k$ and $K$ the metric measure space $(X_{k,\varepsilon},\di_\infty, \mathfrak{m}_{k,K,\varepsilon})$ is a $\CD(0,\infty)$ space for every $\varepsilon<k$. In particular, in the next four sections I will consider the metric measure space $(X_{k,\varepsilon},\di_\infty, \mathfrak{m}_{k,K,\varepsilon})$, but I will avoid to indicate the parameters $k$, $K$ and $\varepsilon$ at the subscript, in order to ease the notation. Later, in the last section, I am going to prove the measured Gromov Hausdorff convergence of the spaces $(X_{k,\varepsilon},\di_\infty, \mathfrak{m}_{k,K,\varepsilon})$ to the space $(X_{k,0},\di_\infty, \mathfrak{m}_{k,K,0})$, as $\varepsilon$ goes to $0$. Combining this with the stability result (Theorem \ref{thm:stabilitycompact}), will follow that $(X_{k,0},\di_\infty, \mathfrak{m}_{k,K,0})$ is itself a $\CD(0,\infty)$ space.

\section{Preliminary Algebraic Lemmas}\label{section:algebraiclemmas}

In this section I simply state and prove four algebraic lemmas that will be fundamental in the proof of the main theorem. In particular they will only help in carrying on the computation, but they do not hide any particular or sophisticated idea. For this reason I do not spend a lot of words on them and I immediately go through the proofs.

\begin{lemma}\label{lem:est2}
There exists a constant $C$ such that for every $A\in\setR^+$ and every $\delta$ with $|\delta|<\frac{1}{32}$ it holds 
\begin{equation} \label{eq:logineq}
    \frac12 \log(A) \leq \log \bigg(1+\bigg(\frac12+\delta\bigg)(A-1)\bigg) + C \delta^2
\end{equation}
\end{lemma}

\begin{pr}
Notice that for some $A\in\setR^+$ the inequality \eqref{eq:logineq} holds also without the term $+C\delta^2$. In particular this is true if and only if 
\begin{equation}\label{eq:1}
    \sqrt{A}\leq 1+\bigg(\frac12+\delta\bigg)(A-1).
\end{equation}
Through elementary computation it is easy to show that inequality \eqref{eq:1} holds for every suitable $\delta$ if
\begin{equation*}
    A \notin \bigg[1-\frac{8|\delta|}{(2\delta+1)^2},1+\frac{8|\delta|}{(2\delta+1)^2} \bigg] \subset [1-16|\delta|,1+16|\delta|]. 
\end{equation*}
Therefore, in order to conclude the proof, it is sufficient to prove inequality \eqref{eq:logineq} for $A\in [1-16|\delta|,1+16|\delta|]$. In this case
\begin{align*}
    \log \bigg(1+\bigg(\frac12+\delta\bigg)(A-1)\bigg) &= \log \bigg(1+\frac12(A-1) + \delta (A-1)\bigg) \\
    & = \log \bigg(1+\frac12(A-1) \bigg) + \int_0^{\delta(A-1)} \frac{1}{1+\frac{1}{2}(A-1)+t} \de t \\
    & \geq \log \bigg(1+\frac12(A-1) \bigg) - \int_0^{|\delta(A-1)|} \frac{1}{1-8|\delta|- 16 \delta^2} \de t \\
    &\geq \log \bigg(1+\frac12(A-1) \bigg) - 2|\delta(A-1)|\\
    &\geq \frac12 \log(A)- 32 \delta^2,
\end{align*}
where the last inequality follows from the concavity of the logarithm.  
\end{pr}

\begin{lemma}\label{prop:basicinterpolats}
Given $(x_0,y_0),(x_1,y_1)\in X$, let $\gamma:[0,1]\to X$ be the function:
\begin{equation*}
    t \mapsto \bigg( (1-t)\frac{y_0}{f(x_0)}  + t \frac{y_1}{f(x_1)}  \bigg) f((1-t)x_0+ t x_1)
\end{equation*}
then
\begin{enumerate}
    \item for every $t\in[0,1]$ it holds
    \begin{equation*}
    \modu{\frac{\gamma'(t)}{x_1-x_0}- \frac{y_1-y_0}{x_1-x_0}}\leq 3k,
\end{equation*}
    \item if $k<\frac14$, for every $t\in[0,1]$ it holds
    \begin{equation*}
    \modu{\gamma''(t)} \leq (x_1-x_0)^2 \frac{2k}{f((1-t)x_0 +t x_1)} \left( \modu{\frac{y_1-y_0}{x_1-x_0}} +2  \right).
\end{equation*}
\end{enumerate}
\end{lemma}

\begin{pr}
The first derivative of $\gamma$ is 
\begin{equation*}
    \gamma'(t)= \bigg( \frac{y_1}{f(x_1)} - \frac{y_0}{f(x_0)} \bigg) f((1-t)x_0+ t x_1)+ \bigg( (1-t)\frac{y_0}{f(x_0)}  + t \frac{y_1}{f(x_1)}  \bigg) f'((1-t)x_0+ t x_1)(x_1-x_0).
\end{equation*}
Therefore, denoting $x_t=(1-t)x_0+tx_1$, it holds 
\begin{align*}
    \frac{\gamma'(t)}{x_1-x_0}&= \frac{y_1-y_0}{x_1-x_0} + \bigg( \frac{1}{f(x_1)}-\frac{1}{f(x_t)} \bigg) \frac{y_1}{x_1-x_0}f(x_t) +  \bigg( \frac{1}{f(x_t)}-\frac{1}{f(x_0)} \bigg) \frac{y_0}{x_1-x_0} f(x_t)\\
    &\quad +\bigg( (1-t)\frac{y_0}{f(x_0)}  + t \frac{y_1}{f(x_1)}  \bigg) f'(x_t)\\
    &= \frac{y_1-y_0}{x_1-x_0} + \frac{y_1}{f(x_1)}\frac{f(x_t)-f(x_1)}{x_1-x_0}+ \frac{y_0}{f(x_0)}\frac{f(x_0)-f(x_t)}{x_1-x_0} \\
    &\quad +\bigg( (1-t)\frac{y_0}{f(x_0)}  + t \frac{y_1}{f(x_1)}  \bigg) f'(x_t).
\end{align*}
On the other hand, it is possible to perform the following estimate 
\begin{align*}
    \modu{\frac{y_1}{f(x_1)}\frac{f(x_t)-f(x_1)}{x_1-x_0}} &\leq\modu{ \frac{y_1}{f(x_1)}}\modu{\frac{f(x_t)-f(x_1)}{x_1-x_0}}\leq \modu{\frac{f(x_t)-f(x_1)}{x_1-x_0}} \leq \sup f'\leq k,
\end{align*}
and the same calculation can be done for the symmetric term, thus
\begin{equation*}
    \modu{\frac{y_0}{f(x_0)}\frac{f(x_0)-f(x_t)}{x_1-x_0}}\leq k.
\end{equation*}
Moreover, a similar procedure shows that
\begin{equation*}
    \modu{\bigg( (1-t)\frac{y_0}{f(x_0)}  + t \frac{y_1}{f(x_1)}  \bigg) f'(x_t)}\leq k.
\end{equation*}
Putting together all this estimates, it is possible to conclude that 
\begin{equation*}
    \modu{\frac{\gamma'(t)}{x_1-x_0}- \frac{y_1-y_0}{x_1-x_0}}\leq 3k.
\end{equation*}
The second derivative of $\gamma$ is 
\begin{align*}
    \gamma''(t)= 2 \bigg( \frac{y_1}{f(x_1)} - \frac{y_0}{f(x_0)} \bigg) f'(x_t)(x_1-x_0) + \bigg( (1-t)\frac{y_0}{f(x_0)}  + t \frac{y_1}{f(x_1)}  \bigg)f''(x_t)(x_1-x_0)^2
\end{align*}
Consider the first term of the right hand side, through algebraic manipulations similar to the ones performed in the first part of the proof, I obtain
\begin{align*}
    2 \bigg( \frac{y_1}{f(x_1)} - &\frac{y_0}{f(x_0)} \bigg) f'(x_t)(x_1-x_0) \\
    &=2 \frac{f'(x_t)(x_1-x_0)^2}{f(x_t)}\bigg[\frac{y_1-y_0}{x_1-x_0} + \frac{y_1}{f(x_1)}\frac{f(x_t)-f(x_1)}{x_1-x_0}+ \frac{y_0}{f(x_0)}\frac{f(x_0)-f(x_t)}{x_1-x_0}\bigg].
\end{align*}
Using the same estimates as before, I get
\begin{align*}
    \modu{2 \bigg( \frac{y_1}{f(x_1)} - \frac{y_0}{f(x_0)} \bigg) f'(x_t)(x_1-x_0)} &\leq 2 (x_1-x_0)^2 \modu{\frac{f'(x_t)}{f(x_t)}} \bigg( \modu{\frac{y_1-y_0}{x_1-x_0}}+ 2k\bigg)\\
    & \leq \frac{2k(x_1-x_0)^2}{f(x_t)} \bigg( \modu{\frac{y_1-y_0}{x_1-x_0}}+ 2k\bigg).
\end{align*}
On the other hand, also the second term of the right hand side can be easily estimate:
\begin{equation*}
    \modu{\bigg( (1-t)\frac{y_0}{f(x_0)}  + t \frac{y_1}{f(x_1)}  \bigg)f''(x_t)(x_1-x_0)^2} \leq (x_1-x_0)^2.
\end{equation*}
Adding this last two inequalities, and using that $f(x_t)\leq 3 k$ I can conclude
\begin{equation*}
    \modu{\gamma''(t)} \leq (x_1-x_0)^2 \bigg[ 1+ \frac{2k}{f(x_t)} \bigg( \modu{\frac{y_1-y_0}{x_1-x_0}}+ 2k\bigg) \bigg] \leq (x_1-x_0)^2 \frac{2k}{f(x_t)} \left( \modu{\frac{y_1-y_0}{x_1-x_0}} +2  \right) ,
\end{equation*}
where the last inequality holds if $k<\frac14$.
\end{pr}

\begin{lemma}\label{lem:est3}
Given a fixed a constant $H$, let $y:I=[x_0,x_1] \to \setR^+$ be a $C^2$ function such that $y'(x)\geq \frac14$ and $y''(x)\leq H \frac{k}{f(x)}$ for every $x \in I$. Then, for $k$ small enough, it holds 
\begin{equation*}
    \log\left(m\left(\frac{x_0+x_1}{2}, y \left(\frac{x_0+x_1}{2} \right) \right)\right) \geq \log(m(x_0,y(x_0))) + \log(m(x_1,y(x_1))) + \frac{K}{128 f(x_1)^2}(x_1-x_0)^2 
\end{equation*}
\end{lemma}

\begin{pr}
Before going into the proof, I want to point out that the inequality I have to prove is basically a $K$-convexity inequality. The strategy of the proof consists in deducing the $K$-convexity from a second derivative estimate. So, start take the first derivative:
\begin{align*}
    \frac{\partial}{\partial x} \left( \log(f(x)) +K\left(\frac{y(x)}{f(x)} \right)^2 \right) = \frac{f'(x)}{f(x)} +2K \frac{y(x)}{f(x)} \left( \frac{y'(x)}{f(x)} - \frac{y(x)f'(x)}{f(x)^2} \right)
\end{align*}
Taking another derivative, I obtain
\begin{align*}
    \frac{\partial^2}{\partial x^2} \left( \log(f(x)) +K\left(\frac{y(x)}{f(x)} \right)^2 \right) &= \frac{f''(x)}{f(x)}-\frac{f'(x)^2}{f(x)^2} +2K \left( \frac{y'(x)}{f(x)} - \frac{y(x)f'(x)}{f(x)^2} \right)^2 \\
    &\, +2K \frac{y(x)}{f(x)} \left( \frac{y''(x)}{f(x)} - 2\frac{y'(x)f'(x)}{f(x)^2} - \frac{y(x)f''(x)}{f(x)^2} + 2 \frac{y(x)f'(x)^2}{f(x)^3} \right).
\end{align*}
Therefore, neglecting some positive terms, the following estimate holds
\begin{align*}
    \frac{\partial^2}{\partial x^2} \left( \log(f(x)) +K\left(\frac{y(x)}{f(x)} \right)^2 \right) &\geq 2K \frac{y'(x)^2}{f(x)^2} -\frac{f'(x)^2}{f(x)^2}-8K \modu{\frac{y(x)y'(x)f'(x)}{f(x)^3}} -2K \modu{\frac{y(x)y''(x)}{f(x)^2}}\\
    & -2K \modu{\frac{y(x)^2f''(x)}{f(x)^3}}.
\end{align*}
Noticing that $\modu{\frac{y(x)}{f(x)}}\leq 1$, $\modu{f'(x)}\leq k$ and $\modu{f(x)}\leq 3k$, I deduce
\begin{align*}
    \frac{\partial^2}{\partial x^2} &\left( \log(f(x)) +K\left(\frac{y(x)}{f(x)} \right)^2 \right)  \\
    &\qquad\qquad\qquad\qquad\geq \frac{2K}{f(x)^2} \bigg( y'(x)^2-\frac{k^2}{2K}-4k\modu{y'(x)}-\modu{f(x)}\modu{y''(x)}-\modu{f(x)}\modu{f''(x)}\bigg)\\
    &\qquad\qquad\qquad\qquad \geq \frac{2K}{f(x)^2} \bigg( y'(x)^2-k^2-4k\modu{y'(x)}-kH-3k\bigg)\\
    &\qquad\qquad\qquad\qquad\geq \frac{K}{16 f(x)^2},
\end{align*}
where the last inequality holds for every suitably small $k$. The thesis follows by making the uniform convexity explicit and noticing that $f(x_1)\geq f(x)$ for every $x\in I$. 
\end{pr}

\noindent Performing the same computations of the previous proof, it possible to prove the following corollary.

\begin{corollary}\label{cor:Kconve}
Given a fixed a constant $H$, let $y:I=[x_0,x_1] \to \setR^+$ be a $C^2$ function such that $y'(x)\geq \frac14$ and $y''(x)\leq H \frac{k}{f(x)}$ for every $x \in I$. Then, for $k$ small enough, it holds 
\begin{equation*}
    K \left(\frac{y \left(\frac{x_0+x_1}{2} \right)}{f \left(\frac{x_0+x_1}{2}\right)}\right)^2 \leq \frac K2 \left( \frac{y(x_0)}{f(x_0)} \right)^2+ \frac K2 \left( \frac{y(x_1)}{f(x_1)} \right)^2 - \frac{K}{128 f(x_1)^2}(x_1-x_0)^2 .
\end{equation*}
\end{corollary}

\begin{lemma}\label{lem:est1}
Let $x\in [-1,1]$ and $\delta>0$ such that $x-\delta,x+\delta \in [-1,1]$, then
\begin{equation*}
    \modu{ \frac{f(x)}{f(x-\delta)} + \frac{f(x)}{f(x+\delta)}-2} \leq \frac{[2k^2 + f(x)]\delta^2}{f(x-\delta)f(x+\delta)}.
\end{equation*}
\end{lemma}

\begin{pr}
Denote $I_1= \int_x^{x-\delta}f'(t)\de t$ and $I_2= \int_x^{x+\delta}f'(t)\de t$, then 
\begin{align*}
       \modu{ \frac{f(x)}{f(x-\delta)} + \frac{f(x)}{f(x+\delta)}-2} &= \modu{ \frac{f(x)}{f(x)+I_1} + \frac{f(x)}{f(x)+I_2}-2}\\
       &=\modu{\frac{-2I_1I_2 -f(x)(I_1+I_2)}{(f(x)+I_1)(f(x)+I_2)}}\\
       &\leq \modu{\frac{2I_1I_2 }{(f(x)+I_1)(f(x)+I_2)}}+\modu{\frac{f(x)(I_1+I_2)}{(f(x)+I_1)(f(x)+I_2)}}.
\end{align*}
But the following estimates hold 
\begin{equation*}
    \modu{I_1+I_2}= \modu{\int_0^\delta f'(x+t)-f'(x-\delta+t) \de t} \leq \int_0^\delta\modu{ f'(x+t)-f'(x-\delta+t)} \de t \leq \delta^2 \sup\modu{f''}\leq \delta^2
\end{equation*}
and 
\begin{equation*}
    \modu{I_1}, \modu{I_2} \leq k \delta.
\end{equation*}
Using this last two estimates I conclude 
\begin{equation*}
    \modu{ \frac{f(x)}{f(x-\delta)} + \frac{f(x)}{f(x+\delta)}-2} \leq \frac{[2k^2 + f(x)]\delta^2}{(f(x)+I_1)(f(x)+I_2)}
\end{equation*}
\end{pr}

\section{How to Prove Convexity of the Entropy}\label{sec:jacobi}

In this section I prove an important result (Proposition \ref{prop:jacobi}) that will be fundamental in the following, in order to prove the CD condition. This results relies on the possibility to  compute the density of a pushforward measure, through Jacobi equation. For example, consider two measures $\mu_0,\mu_1\in \ProbTwo(\R^2)$ which are absolutely continuous with respect to the Lebesgue measure $\Leb^2$, with density $\rho_0$ and $\rho_1$ respectively. Suppose there exists a smooth one-to-one map $T:\R^2 \to \R^2$ such that $T_\#\mu_0=\mu_1$, then it is well known that 
\begin{equation} \label{eq:jacobi}
    \rho_1(T(x,y)) J_T(x,y)= \rho_0(x,y),
\end{equation}
for $\mu_0$-almost every $(x,y)$. As shown in \cite{ambrosio2005gradient} the assumptions on the map $T$ can be relaxed and it is sufficient to require $T$ to be approximately differentiable and injective outside a $\mu_0$-null set. However, in this work I am going to deal with transport map which are not necessarily approximately differentiable, but have another rigidity property. Therefore, the version of Jacobi equation I will need is the following, which is an easy consequence of Proposition 2.1 in \cite{rajala2013failure}.

\begin{prop}\label{prop:rajalajacobi}
Let $\mu_0,\mu_1\in \ProbTwo(\R^2)$ be absolutely continuous with respect to the Lebesgue measure $\Leb^2$. Assume that there exists a map $T=(T_1,T_2)$ which is injective outside a $\mu_0$-null set, such that $T_\# \mu_0=\mu_1$. Suppose also that $T_1$ locally does not depend on the $y$ coordinate and it is increasing in $x$, while $T(x,y)$ is increasing in $y$ for every fixed $x$. Then the Jacobi equation \eqref{eq:jacobi} is satisfied with $J_T=\frac{\partial T_1}{\partial x} \frac{\partial T_2}{\partial y}$.
\end{prop}

Now that I have shown that Jacobi equation can be adapted to the context of this work, I am going to explain how it can be useful in proving convexity of the entropy functional. The following proposition will be an important element in the main proof of this article. 

\begin{prop} \label{prop:jacobi}
Let $\mu_0,\mu_1\in \ProbTwo(X)$ be absolutely continuous measure and let $T:X \to X$ be an optimal transport map between $\mu_0$ and $\mu_1$, in particular $T_\# \mu_0 = \mu_1$. Consider a midpoint $\mu_{1/2}$ of $\mu_0$ and $\mu_1$, assume that $\mu_{1/2}=[M\circ(\id,T)]_\# \mu_0$ where the map $M:X\times X \to X$ is a midpoint selection. Suppose also that the maps $T$ and $M\circ(\id,T):X \to X$ are injective outside a set of measure zero and they satisfy the Jacobian equation, with suitable Jacobian functions $J_T$ and $J_{M\circ(\id,T)}$. If 
\begin{equation*}
    \log \left( m\big(M((x,y),T(x,y))\big) J_{M\circ(\id,T)}(x,y) \right) \geq \frac12 \log \left( m(T(x,y)) J_T(x,y) \right) + \frac12 \log (m(x,y))
\end{equation*}
for $\mu_0$ almost every $(x,y)$, then
\begin{equation}\label{eq:entconvexity}
    \Ent(\mu_{1/2}) \leq \frac12 \Ent(\mu_0) + \frac12 \Ent(\mu_1).
\end{equation}
\end{prop}

\begin{pr}

Suppose $\mu_0=\rho_0 \mathfrak{m}= \tilde{\rho}_0 \Leb^2$, $\mu_1=\rho_1 \mathfrak{m}= \tilde{\rho}_1 \Leb^2$ and $\mu_{1/2}=\rho_{} \mathfrak{m}= \tilde{\rho}_1 \Leb^2$. It easy to realize that, in order to prove \eqref{eq:entconvexity}, it is sufficient to prove that 
\begin{equation}\label{eq:sufcond}
    \log\big( \rho_{1/2}\big(M((x,y),T(x,y))\big) \big) \leq \frac12 \log \big( \rho_{1}(T(x,y)) \big) + \frac12 \log(\rho_0(x,y)),
\end{equation}
for $\mu_0$-almost every $(x,y)$.
On the other hand, the validity of Jacobi equation ensures that 
\begin{equation*}
    \tilde{\rho}_1(T(x,y)) J_T(x,y) = \tilde{\rho}_0(x,y),
\end{equation*}
for $\mu_0$-almost every $(x,y)$, and thus that 
\begin{equation*}
    \rho_1(T(x,y)) m(T(x,y)) J_T(x,y) = \rho_0(x,y) m(x,y).
\end{equation*}
for $\mu_0$-almost every $(x,y)$.
Analogously, it is possible to infer that
\begin{equation*}
    \rho_{1/2}\big(M((x,y),T(x,y))\big) m\big(M((x,y),T(x,y))\big) J_{M\circ(\id,T)}(x,y) = \rho_0(x,y) m(x,y).
\end{equation*}
for $\mu_0$-almost every $(x,y)$.
Therefore \eqref{eq:sufcond} is equivalent to 
\begin{equation*}
    \log \left( \frac{\rho_0(x,y) m(x,y)}{ m\big(M((x,y),T(x,y))\big) J_{M\circ(\id,T)}(x,y)}\right) \leq 
    \frac12 \log \left( \frac{\rho_0(x,y) m(x,y)}{m(T(x,y)) J_T(x,y)} \right) + \frac12 \log(\rho_0(x,y)).
\end{equation*}
Some easy algebraic computations show that this last equation is equivalent to
\begin{equation*}
    \log \left( m\big(M((x,y),T(x,y))\big) J_{M\circ(\id,T)}(x,y) \right) \geq \frac12 \log \left( m(T(x,y)) J_T(x,y) \right) + \frac12 \log (m(x,y)),
\end{equation*}
concluding the proof.
\end{pr}

\noindent Notice that the result of this last proposition gains importance if seen in relation with Proposition \ref{prop:midpoint}. In fact, Proposition \ref{prop:jacobi} provides a strategy to prove entropy convexity in a suitable midpoint, which is sufficient to deduce CD condition, according to Proposition \ref{prop:midpoint}.

I conclude the section with a simple corollary of Proposition \ref{prop:jacobi}, which will be useful in the last section of this work. This result is a straightforward consequence of the previous proof, and takes full advantage of the fact that Jacobi equation allows to prove entropy convexity pointwise as well as globally.

\begin{corollary}\label{cor:jacobi}
Under the same assumptions of Proposition \ref{prop:jacobi}, it holds that
\begin{equation*}
    \Ent\big([M\circ(\id,T)]_\# (f \mu_0) \big) \leq \frac12 \Ent(f \mu_0) + \frac12 \Ent(T_\#(f\mu_0)),
\end{equation*}
for every bounded measurable function $f:X\to \R^+$, with $\int f \de \mu_0=1$.
\end{corollary}

\section{Definition of the Midpoint}\label{sec:defimidpoint}

As previously pointed out, in order to prove CD condition I am going to prove entropy convexity in a suitable midpoint of any pair of absolutely continuous measures. Notice that in an highly branching metric measure space the choice of a midpoint can be done with great freedom. This is because, in general, both the optimal transport map and the geodesic interpolation are not unique, and thus they must be selected in a clever way. In this section, for any pair of absolutely continuous measures, I define a suitable midpoint and in the following section I will show that it satisfies the convexity of the entropy. This midpoint selection is actually quite complicated but it does the job, in particular I believe there is no way to obtain a considerably simpler one.

Before going into the details, I introduce the following subsets of $\R^2 \times \R^2$ that will play an important role in the definition of the midpoint.

\begin{definition}
Define the sets $V,D,H,H_0,H_\frac12,H_1 \subset \R^2 \times \R^2$ as:
\begin{equation*}
    V:=\left\{\left(\left(x_{0}, y_{0}\right),\left(x_{1}, y_{1}\right)\right) \in \mathbb{R}^{2} \times \mathbb{R}^{2}:\left|x_{0}-x_{1}\right|<\left|y_{0}-y_{1}\right|\right\},
\end{equation*}
\begin{equation*}
    D:=\left\{\left(\left(x_{0}, y_{0}\right),\left(x_{1}, y_{1}\right)\right) \in \mathbb{R}^{2} \times \mathbb{R}^{2}:\left|x_{0}-x_{1}\right|=\left|y_{0}-y_{1}\right|\right\},
\end{equation*}
\begin{equation*}
    H:=\left\{\left(\left(x_{0}, y_{0}\right),\left(x_{1}, y_{1}\right)\right) \in \mathbb{R}^{2} \times \mathbb{R}^{2}:\left|x_{0}-x_{1}\right|>\left|y_{0}-y_{1}\right|\right\}= H_0  \cup H_1,
\end{equation*}
where
\begin{equation*}
    H_0:=\left\{\left(\left(x_{0}, y_{0}\right),\left(x_{1}, y_{1}\right)\right) \in \mathbb{R}^{2} \times \mathbb{R}^{2}:\frac12\left|x_{0}-x_{1}\right|\geq\left|y_{0}-y_{1}\right|\right\},
\end{equation*}
\begin{equation*}
    H_1:=\left\{\left(\left(x_{0}, y_{0}\right),\left(x_{1}, y_{1}\right)\right) \in \mathbb{R}^{2} \times \mathbb{R}^{2}:\left|x_{0}-x_{1}\right|>\left|y_{0}-y_{1}\right| >\frac12 |x_0-x_1|\right\}.
\end{equation*}
\end{definition}

The first step in the selection of the midpoint consists in choosing a suitable optimal transport map between two given absolutely continuous measures. To this aim I use a nice result by Rajala \cite{rajala2013failure}, who was able to show the existence of an optimal transport map with different nice properties. The idea behind his work is to use consecutive minimization, in order to select a particular optimal transport plan. The result of Rajala can be summarized in the following statement.

\begin{prop}\label{prop:map}
Given two measures $\mu_0,\mu_1 \in \Prob(\R^2)$ which are absolutely continuous with respect to the Lebesgue measure $\Leb^2$, there exists a measurable optimal transport map $T=(T_1,T_2)$, injective outside a $\mu_0$-null set, satisfying $T_\#\mu_0 = \mu_1$, with some nice rigidity properties. In particular the optimal transport plan $(\id,T)_\# \mu_0$ is concentrated in a set $\Gamma\subset X \times X$, such that for all $\left(x, y_{1}\right),\left(x, y_{2}\right),\left(x_{1}, y\right),\left(x_{2}, y\right) \in\left\{(x, y) \in X:((x, y), T(x, y)) \in \Gamma\right\}$ it holds that
\begin{equation*}
    \text {if } y_{1} \neq y_{2} \text { and } T_{1}\left(x, y_{1}\right)=T_{1}\left(x, y_{2}\right), \text { then } \frac{T_{2}\left(x, y_{1}\right)-T_{2}\left(x, y_{2}\right)}{y_{1}-y_{2}} \geq 0
\end{equation*}
and
\begin{equation*}
    \text {if } x_{1} \neq x_{2} \text { and } T_{2}\left(x_{1}, y\right)=T_{2}\left(x_{2}, y\right), \text { then } \frac{T_{1}\left(x_{1}, y\right)-T_{1}\left(x_{2}, y\right)}{x_{1}-x_{2}} \geq 0.
\end{equation*}
Moreover for $\mu_0$-almost every $(x,y)$ I have that
\begin{equation*}
    \begin{array}{l}T_{1} \text { is locally constant in } y, \text { if }((x, y), T(x, y)) \in H \text { and } \\ T_{2} \text { is locally constant in } x, \text { if }((x, y), T(x, y)) \in V\end{array}.
\end{equation*}
Combining this two properties with some monotonicity properties one can deduce that the function $T_1(x,y)$ is increasing in $x$ for every fixed $y$ and the function $T_2(x,y)$ is increasing in $y$ for every fixed $x$, as a consequence for $\mu_0$-almost every $(x,y)$ it holds
\begin{equation*}
    \frac{\partial T_{1}}{\partial x} \geq 0 \text { and } \frac{\partial T_{2}}{\partial y} \geq 0, \text { if }((x, y), T(x, y)) \in H \cup V.
\end{equation*}
Finally let me point out that, since $(\id,T)_\# \mu_0\in \OptPlans(\mu_0,\mu_1)$ the usual cyclical monotonicity holds, thus I can assume that for every $\left(z_{1}, w_{1}\right),\left(z_{2}, w_{2}\right) \in \Gamma$ it holds
\begin{equation*}
    \di_\infty^2(z_1,w_1)+\di_\infty^2(z_2,w_2)\leq \di_\infty^2(z_1,w_2)+\di_\infty^2(z_2,w_1).
\end{equation*}
\end{prop}

Now fix two measures $\mu_0,\mu_1\in \Prob(X)$ absolutely continuous with respect to the reference measure $\m$, and thus also with respect to the Lebesgue measure $\Leb^2$. Call $T$ the optimal transport map between $\mu_0$ and $\mu_1$, that satisfies the requirements of Proposition \ref{prop:map}. Moreover denote by $\Gamma$ the full $(\id,T)_\# \mu_0$-measure set with all the monotonicity property stated in Proposition \ref{prop:map}. In order to identify a midpoint of $\mu_0$ and $\mu_1$, I need to choose a proper midpoint interpolation, that is a measurable map $M:X\times X\to X$ such that 
\begin{equation}\label{eq:midpointinter}
    \di_\infty(M(z,w),z) = \di_\infty(M(z,w),w)= \frac12 \di_\infty(z,w) \quad \text{for every }(z,w)\in X \times X ,
\end{equation}
the desired midpoint will be $M_\# \big( (\id,T)_\# \mu_0 \big)= [M\circ (\id,T)]_\# \mu_0$.

Let me now define the midpoint interpolation map $M$ that I will use for the proof of the main theorem. This definition is actually quite involved, in fact the map $M$ is defined in different ways on the sets $V$, $D$, $H_0$ and $H_1$. In particular the precise definition is the following:
\begin{itemize}
    \item If $\big((x_0,y_0), (x_1,y_1)\big)\in V\cup D$
    \begin{equation*}
        M\big((x_0,y_0), (x_1,y_1)\big):=\left(\frac{x_0+x_1}{2}, \frac{y_0+y_1}{2} \right)
    \end{equation*}
    \item If $\big((x_0,y_0), (x_1,y_1)\big)\in H_0 $,
    \begin{equation*}
    M\big((x_0,y_0), (x_1,y_1)\big)= \left(\frac{x_0+x_1}{2} , \frac12 \bigg( \frac{y_0}{f(x_0)}  +  \frac{y_1}{f(x_1)}  \bigg)f\left(\frac{x_0+x_1}{2}\right) \right).\end{equation*}
    \item If $\big((x_0,y_0), (x_1,y_1)\big)\in H_1$, with $x_0<x_1$ and $y_0<y_1$, introduce the quantity
    \begin{equation*}
        \tilde{y}(x_0,x_1,y_0)= \frac12\bigg( \frac{y_0}{f(x)}  +  \frac{y_0+\frac{x_1-x_0}{2}}{f(x_1)}  \bigg)f\left(\frac{x_1+x_0}{2}\right)  -y_0,
    \end{equation*}
    and consequently define 
    \begin{align*}
        M&\big((x_0,y_0), (x_1,y_1)\big)\\
        & \qquad\qquad= \left(\frac{x_0+x_1}{2} ,y_0 + \tilde{y}(x_0,x_1,y_0) + \left(\frac{x_1-x_0}{2} - \tilde{y}(x_0,x_1,y_0) \right) \left( 2 \frac{y_1 - y_0}{x_1-x_0} -1 \right)\right).
    \end{align*}
    In the other cases $M$ can be defined analogously, but I prefer not to explicitly do it, in order to avoid unnecessary complications in this definition and in the following. In particular every proof from now on will be done using only the definition above, implying it can be easily adapted to the other cases.
\end{itemize}
With such a complex definition is not completely obvious that the map $M$ satisfies condition \eqref{eq:midpointinter}, but this can be easily proven.
\begin{prop}
The map $M$ actually defines a midpoint interpolation. 
\end{prop}

\begin{pr}
This is completely obvious if $\big((x_0,y_0), (x_1,y_1)\big)\in V\cup D$. While when $\big((x_0,y_0), (x_1,y_1)\big)\in H_0$ the statement is an easy consequence of Lemma \ref{prop:basicinterpolats}, provided that $k$ is sufficiently small, and of the fact that
\begin{equation*}
    0 \leq \frac12 \bigg( \frac{y_0}{f(x_0)}  +  \frac{y_1}{f(x_1)}  \bigg)f\left(\frac{x_0+x_1}{2}\right) \leq f\left(\frac{x_0+x_1}{2}\right).
\end{equation*}
On the other hand if $\big((x_0,y_0), (x_1,y_1)\big)\in H_1$ it is sufficient to notice that the point 
\begin{equation*}
  (\bar x, \bar y):=\left(\frac{x_0+x_1}{2} ,y_0 + \tilde{y}(x_0,x_1,y_0) + \left(\frac{x_1-x_0}{2} - \tilde{y}(x_0,x_1,y_0) \right) \left( 2 \frac{y_1 - y_0}{x_1-x_0} -1 \right)\right)
\end{equation*}
lies on a convex combination of a 45 degree line and a curve $\gamma:[x_0,x_1]\to X$, with 
\begin{equation*}
    \frac 12 -3k \leq \gamma'(x)\leq \frac 12 +3k \quad \text{ for every } x\in [x_0,x_1],
\end{equation*}
according to Lemma \ref{prop:basicinterpolats}. As a consequence, for a suitably small $k$, this ensures both that $(\bar x, \bar y)\in X$ and that $(\bar x, \bar y)$ is a $\di_\infty$-midpoint of $(x_0,y_0)$ and $(x_1,y_1)$.
\end{pr}

In order to efficiently apply Proposition \ref{prop:jacobi}, I need to compute the Jacobian of the map $M\circ (\id,T)$. Observe that the way $M$ is defined, combined with the properties of $T$, ensures that $M\circ (\id,T)$ satisfies all the rigidity assumptions of Proposition \ref{prop:rajalajacobi}. Therefore, proving the following result, will allow to use Proposition \ref{prop:rajalajacobi} to compute $J_{M\circ (\id,T)}$.

\begin{prop}\label{prop:aeinjectivity}
The map $M\circ(\id,T)$ is injective outside a $\mu_0$-null set.
\end{prop}

\begin{pr}
First of all notice that it is sufficient to prove the injectivity of $M$ on $\Gamma$, because $\Gamma$ has full $(\id,T)_\# \mu_0$-measure. Thus suppose by contradiction that there exist $$((x_{1}, y_{1}), (x_{1}^{\prime}, y_{1}^{\prime}))\ne((x_{2}, y_{2}),(x_{2}^{\prime}, y_{2}^{\prime})) \in \Gamma$$ such that \begin{equation*}
    M((x_{1}, y_{1}), (x_{1}^{\prime}, y_{1}^{\prime}))= M((x_{2}, y_{2}),(x_{2}^{\prime}, y_{2}^{\prime})).
\end{equation*}
Following the proof of Lemma 3.7 in \cite{rajala2013failure}, I can limit myself to consider the case when $$((x_{1}, y_{1}), (x_{1}^{\prime}, y_{1}^{\prime}))\ne((x_{2}, y_{2}),(x_{2}^{\prime}, y_{2}^{\prime})) \in \Gamma \cap H.$$ In this case, one can easily realize that the cyclical monotonicity properties of $\Gamma$ imply that $x_1=x_2$ and $x_1'=x_2'$. So, if $((x_{1}, y_{1}), (x_{1}^{\prime}, y_{1}^{\prime})),((x_{2}, y_{2}),(x_{2}^{\prime}, y_{2}^{\prime})) \in \Gamma \cap H_0$, thesis simply follows from the definition of $M$ and from the monotonicity of $T_2$. While if $((x_{1}, y_{1}), (x_{1}^{\prime}, y_{1}^{\prime})),((x_{2}, y_{2}),(x_{2}^{\prime}, y_{2}^{\prime})) \in \Gamma \cap H_1$ the statement is a consequence of the monotonicity property of $T_2$ associated with the fact that the quantity 
\begin{equation*}
    y + \tilde{y}(x,x',y) + \left(\frac{x'-x}{2} - \tilde{y}(x,x',y) \right) \left( 2 \frac{y' - y}{x'-x} -1 \right)
\end{equation*}
is locally increasing in $y$ and $y'$, when $((x,y),(x',y'))\in H_1$ (with $x<x'$ and $y<y'$). The first monotonicity is not straightforward, therefore I am going to prove it. First of all notice that, according to Lemma \ref{lem:est1} and since $f\leq 3k$, it holds that \begin{equation*}
    2 \modu{\frac{\partial}{\partial y}\tilde{y}(x,x',y)} = \modu{ \frac{f\big(\frac{x+x'}{2}\big)}{f(x)} + \frac{f\big(\frac{x+x'}{2}\big)}{f(x')}-2} \leq \frac{[2k^2 + 3k]\big(\frac{x'-x}{2}\big)^2}{f(x)f(x')}.
\end{equation*}
Moreover, the geometry of the set $X$ allows to deduce that
\begin{equation*}
    \frac{x'-x}{2} \leq y + \frac{x'-x}{2} \leq f(x') \leq f(x) + k (x'-x),
\end{equation*}
and consequently $\big( 1 -2k\big) \frac{x'-x}{2} \leq f(x)\leq f(x')$. On the other hand, observe that Lemma \ref{prop:basicinterpolats} guarantees that
\begin{equation*}
    \left(\frac{x'-x}{2} - \tilde{y}(x,x',y) \right) \leq \left( \frac 12 + 3k  \right)\frac{x'-x}{2}
\end{equation*}
Those estimates allow to conclude that 
\begin{align*}
    \frac{\partial}{\partial y}  \bigg[ y + \tilde{y}(x,x',y) &+ \left(\frac{x'-x}{2} - \tilde{y}(x,x',y) \right) \left( 2 \frac{y' - y}{x'-x} -1 \right) \bigg] \\
    &\geq \frac 12 - 3k - \modu{\frac{\partial}{\partial y}\tilde{y}(x,x',y)} \left( 2 - 2 \frac{y' - y}{x'-x}  \right) \\
    & \geq \frac 12 - 3k - \frac{[2k^2 + 3k]}{(1-2k)^2} >0,
\end{align*}
for $k$ sufficiently small. 
The case when $((x_{1}, y_{1}), (x_{1}^{\prime}, y_{1}^{\prime}))\in \Gamma \cap H_1,((x_{2}, y_{2})$ and $(x_{2}^{\prime}, y_{2}^{\prime})) \in \Gamma \cap H_0$ can be treated analogously.
\end{pr}

%EDIT cambiare la stima su f\leq k

\section{Proof of CD Condition}\label{section:CDproof}

In the previous sections I have introduced all I need to prove that the metric measure space $(X,\di_\infty,\m)$ satisfies the $\CD(0,\infty)$ condition. Let me now go into the details of the proof.

\begin{theorem}\label{thm:CDTheExample}
For suitable $k$ and $K$, the metric measure space $(X_{k,\varepsilon},\di_\infty, \mathfrak{m}_{k,K,\varepsilon})$ is a $\CD(0,\infty)$ space, for every $0<\varepsilon<k$.
\end{theorem}

\begin{pr}
Let $\mu_0,\mu_1\in \Prob(X)$ be absolutely continuous with respect to the reference measure $\m$, then, according to Proposition \ref{prop:midpoint}, it is sufficient to prove that 
\begin{equation*}
    \Ent\big([M\circ (\id,T)]_\# \mu_0\big) \leq \frac12 \Ent(\mu_0) + \frac12 \Ent(\mu_1).
\end{equation*}
Given Proposition \ref{prop:jacobi}, it is enough check the validity of 
\begin{equation}\label{eq:condition}
    \log \left( m\big(M((x,y),T(x,y))\big) J_{M\circ(\id,T)}(x,y) \right) \geq \frac12 \log \left( m(T(x,y)) J_T(x,y) \right) + \frac12 \log (m(x,y))
\end{equation}
for $\mu_0$-almost every $(x,y)$. 
For $\mu_0$-almost every $(x,y)\in V \cup D$ this can be done following \cite{rajala2013failure}. Thus I will treat the other cases and applying Lemma \ref{lem:est3}. \\
Notice that, for $\mu_0$-almost every $(x,y)$ such that $((x,y),T(x,y))\in H_0 \cap \Gamma$, I have
\begin{equation*}
    M\circ(\id,T) (x,y)=  \left(\frac{x+T_1}{2} , \frac12 \bigg( \frac{y}{f(x)}  +  \frac{T_2}{f(T_1)}  \bigg)f\left(\frac{x+T_1}{2}\right) \right).
\end{equation*}
Then, according to what I did in previous sections, it is possible to apply Proposition \ref{prop:rajalajacobi} and deduce that
\begin{equation*}
    J_{M\circ(\id,T)} (x,y) = \frac{1}{2} \left( 1 + \frac{\partial T_1}{\partial x} \right) \frac{1}{2} \bigg( \frac{1}{f(x)}  +  \frac{\frac{\partial T_2}{\partial y}}{f(T_1)}  \bigg) f\left( \frac{x+T_1}{2} \right),
\end{equation*}
for $\mu_0$-almost every $(x,y)$ such that $((x,y),T(x,y))\in H_0 \cap \Gamma$.
Furthermore it holds that
\begin{equation*}
    m\big(M((x,y),T(x,y))\big) = f\left( \frac{x+T_1}{2} \right)^{-1} \exp\left( \frac{-K}{4} \bigg( \frac{y}{f(x)}  +  \frac{T_2}{f(T_1)}  \bigg)^2\right),
\end{equation*}
thus, putting together this last two relations, I obtain
\begin{align*}
    \log \left( m\big(M((x,y),T(x,y))\big) J_{M\circ(\id,T)}(x,y) \right) &= \log \Bigg( \frac{1}{2} \left( 1 + \frac{\partial T_1}{\partial x} \right) \frac{1}{2} \bigg( \frac{1}{f(x)}  +  \frac{\frac{\partial T_2}{\partial y}}{f(T_1)}  \bigg) \\ 
    &\qquad\exp\Bigg( \frac{-K}{4} \bigg( \frac{y}{f(x)}  +  \frac{T_2}{f(T_1)}  \bigg)^2\Bigg) \Bigg)\\
    &= \log \left(\frac{1}{2} \left( 1 + \frac{\partial T_1}{\partial x} \right) \right) + \log \left( \frac{1}{2} \bigg( \frac{1}{f(x)}  +  \frac{\frac{\partial T_2}{\partial y}}{f(T_1)}  \bigg)\right)\\
    &\quad - K  \bigg(\frac12 \bigg( \frac{y}{f(x)}  +  \frac{T_2}{f(T_1)}  \bigg) \bigg)^2
\end{align*}
On the other hand it holds 
\begin{equation*}
    \log(m(x,y)) = \log\left(\frac{1}{f(x)} \exp \left(-K \bigg(\frac{y}{f(x)}\bigg)^2 \right) \right)= \log(1)+ \log\left(\frac{1}{f(x)} \right) -K \bigg(\frac{y}{f(x)}\bigg)^2 
\end{equation*}
and, applying once again Proposition \ref{prop:rajalajacobi}, this time to the map $T$, also
\begin{align*}
    \log \left( m(T(x,y)) J_T(x,y) \right)&= \log \left( \frac{\partial T_1}{\partial x} \frac{\partial T_2}{\partial y} \frac{1}{f(T_1)} \exp \left(-K \bigg(\frac{T_2}{f(T_1)}\bigg)^2 \right)\right)\\
    &=\log \left( \frac{\partial T_1}{\partial x} \right)+ \log \left( \frac{\frac{\partial T_2}{\partial y}}{f(T_1)}\right)-K \bigg(\frac{T_2}{f(T_1)}\bigg)^2,
\end{align*}
for $\mu_0$-almost every $(x,y)$ such that $((x,y),T(x,y))\in H_0 \cap \Gamma$.
Putting together this last three equations, inequality \eqref{eq:condition} follows from the concavity of the functions $\log$ and $-Kx^2$.
\\
Passing now to the last case, for $\mu_0$-almost every $(x,y)$ such that $((x,y),T(x,y))\in H_1 \cap \Gamma$ (with $x<T_1(x,y)$ and $y<T_2(x,y)$) I have 
\begin{align*}
    (S_1,S_2)(x,y)&:= M \circ (\id, T)(x,y) \\
   &\quad = \left(\frac{x+T_1}{2}, y + \tilde{y}(x,T_1,y) + \left(\frac{T_1-x}{2} - \tilde{y}(x,T_1,y) \right) \left( 2 \frac{T_2 - y}{T_1-x} -1 \right) \right).
\end{align*}
Reasoning as before, Proposition \ref{prop:rajalajacobi} ensures that, for $\mu_0$-almost every $(x,y)$ such that $((x,y),T(x,y))\in H_1 \cap \Gamma$ (with $x<T_1(x,y)$ and $y<T_2(x,y)$),
\begin{equation*}
     J_{M\circ(\id,T)} (x,y)= \frac{\partial S_1}{\partial x} \frac{\partial S_2}{\partial y},
\end{equation*}
and in particular it holds that   
\begin{equation*}
    \frac{\partial S_1}{\partial x}=\frac{1}{2} \left( 1 + \frac{\partial T_1}{\partial x} \right),
\end{equation*}
and 
\begin{align*}
    \frac{\partial S_2}{\partial y} &= 1+ \frac{\partial}{\partial y}\tilde{y}(x,T_1,y) \left( 2 - 2 \frac{T_2 - y}{T_1-x} \right) + 2\left(\frac{T_1-x}{2} - \tilde{y}(x,T_1,y) \right) \frac{\frac{\partial T_2}{\partial y}-1}{T_1-x}\\
    &= 1+ \frac{\partial}{\partial y}\tilde{y}(x,T_1,y) \left( 2 - 2 \frac{T_2 - y}{T_1-x} \right) + \left( \frac{\partial T_2}{\partial y}-1 \right) \left(1 - \frac{\tilde{y}(x,T_1,y)}{\frac{T_1-x}{2}} \right) \\
    &= 1+ \frac{\partial}{\partial y}\tilde{y}(x,T_1,y) \left( 2 - 2 \frac{T_2 - y}{T_1-x} \right) + \left( \frac{\partial T_2}{\partial y}-1 \right) \left(\frac12 - \frac{\tilde{y}(x,T_1,y)-\frac{T_1-x}{4}}{\frac{T_1-x}{2}} \right).
\end{align*}
I can now consider the explicit value of $\tilde{y}(x,T_1,y)$ and notice that
\begin{align*}
    \tilde{y}(x,T_1,y)-\frac{x+T_1}{4}&= \frac12 \left[  \bigg( \frac{y}{f(x)}  +  \frac{y+\frac{T_1-x}{2}}{f(T_1)}  \bigg)f\left(\frac{T_1+x}{2}\right)  -2y -\frac{T_1-x}{2} \right]\\
    & = \frac12 \left[ y \left( \frac{f\left(\frac{T_1+x}{2} \right)}{f(x)} + \frac{f\left(\frac{T_1+x}{2} \right)}{f(T_1)}-2 \right) + \frac{T_1-x}{2} \left( \frac{f\left(\frac{T_1+x}{2} \right)}{f(T_1)} -1 \right)\right].
\end{align*}
Moreover, I can easily obtain that
\begin{equation*}
    \modu{\frac{f\left(\frac{T_1+x}{2} \right)}{f(T_1)} -1 } = \modu{\frac{f\left(\frac{T_1+x}{2} \right)-f(T_1)}{f(T_1)} } \leq \sup f' \cdot \frac{\frac{T_1-x}{2}}{f(T_1)} =  k \cdot \frac{\frac{T_1-x}{2}}{f(T_1)} ,
\end{equation*}
thus applying Lemma \ref{lem:est1} and noticing that $\frac{y}{f(x)}\leq1$, I can conclude that
\begin{align*}
    \modu{\tilde{y}(x,T_1,y)-\frac{T_1-x}{4}} \leq \frac{\left[2k^2 + k + f\left(\frac{T_1+x}{2} \right)\right]\left(\frac{T_1-x}{2}\right)^2 }{f(T_1)} \leq\frac{[2k^2 + 4k]\left(\frac{T_1-x}{2}\right)^2}{f(T_1)}
\end{align*}
where in the last inequality I used that $f(x)\leq 3k$ for every $x\in[-1,1]$. Therefore it is possible to infer that 
\begin{equation}\label{eq:firstest}
    \modu{\frac{\tilde{y}(x,T_1,y)-\frac{T_1-x}{4}}{\frac{T_1-x}{2}}} \leq \frac{[2k^2 + 4k]\frac{T_1-x}{2}}{f(T_1)} \leq \frac{1}{64} \frac{\frac{T_1-x}{2}}{f(T_1)}\leq \frac{\frac{T_1-x}{2}}{f(T_1)},
\end{equation}
for a sufficiently small $k$. 
Now suppose that 
\begin{equation*}
    \frac{\partial}{\partial y}\tilde{y}(x,T_1,y)=  \frac 12 \frac{f\big(\frac{x+T_1}{2}\big)}{f(x)} + \frac 12 \frac{f\big(\frac{x+T_1}{2}\big)}{f(T_1)}-1 >0 ,
\end{equation*}
then, after noticing that the geometry of the set $X$ ensures that $\frac{T_1-x}{2}\leq f(T_1)$, it is  possible to apply Lemma \ref{lem:est2} and obtain that
\begin{equation}\label{eq:first}
\begin{aligned}
    \log \left( \frac{\partial S_2}{\partial y} \right) &\geq \log \left( 1+ \left( \frac{\partial T_2}{\partial y}-1 \right) \left(\frac12 - \frac{\tilde{y}(x,T_1,y)-\frac{T_1-x}{4}}{\frac{T_1-x}{2}} \right) \right)\\
    &\geq  \log(1) + \log \left( \frac{\partial T_2}{\partial y} \right)-C \frac{\left(\frac{T_1-x}{2}\right)^2}{f(T_1)^2} .
\end{aligned}
\end{equation}
On the other hand, it is easy to realize that the point $M \circ (\id, T)(x,y)$ lies on a curve, which is a convex combination of a 45 degree line and of the curve 
\begin{equation*}
    t \mapsto \bigg( (1-t)\frac{y}{f(x)}  + t \, \frac{y+ \frac{T_1-x}{2}}{f(T_1)}  \bigg) f((1-t)x+ t T_1).
\end{equation*}
Therefore, up to take a suitably small $k$, Lemma \ref{prop:basicinterpolats} allows to apply Lemma \ref{lem:est3}, and obtain
\begin{equation}\label{eq:second}
    \log\left(m\left(M \circ (\id, T)(x,y)  \right)\right) \geq \log(m(x,y)) + \log(m(T(x,y)))+ \frac{K}{128 f(T_1)^2}(T_1-x)^2.
\end{equation}
Inequality \eqref{eq:condition} then follows as before, putting together \eqref{eq:first} and \eqref{eq:second} and taking $K$ sufficiently large.\\
Suppose instead that 
\begin{equation*}
    \frac{\partial}{\partial y}\tilde{y}(x,T_1,y)=  \frac 12 \frac{f\big(\frac{x+T_1}{2}\big)}{f(x)} + \frac 12 \frac{f\big(\frac{x+T_1}{2}\big)}{f(T_1)}-1 <0 ,
\end{equation*}
then notice that 
\begin{equation}\label{eq:negativecase1}
    \begin{aligned}
         \log \left( f\bigg(\frac{T_1+x}{2}\bigg) \right) - \frac 12 \big( \log (f(x)) + \log (f(T_1))\big) = \frac 12 \log \left( \frac{f\left(\frac{T_1+x}{2} \right)}{f(x)} \cdot \frac{f\left(\frac{T_1+x}{2} \right)}{f(T_1)} \right) \\
        \leq  \log \left( \frac 12 \frac{f\left(\frac{T_1+x}{2} \right)}{f(x)} + \frac 12 \frac{f\left(\frac{T_1+x}{2} \right)}{f(T_1)} \right) =  \log \left( 1 + \frac{\partial}{\partial y}\tilde{y}(x,T_1,y) \right)
    \end{aligned}
\end{equation}
Moreover, according to the estimates done in the proof of Proposition \ref{prop:aeinjectivity} and to \eqref{eq:firstest}, it is easy to realize that, for $k$ small enough
\begin{equation*}
    \left(\frac12 - \frac{\tilde{y}(x,T_1,y)-\frac{T_1-x}{4}}{\frac{T_1-x}{2}} \right)\left( 1 + \frac{\partial}{\partial y}\tilde{y}(x,T_1,y) \right)^{-1} = \frac 12 + \tilde \delta 
\end{equation*}
for some $\tilde \delta$ such that $$\big|\tilde \delta\big| \leq \frac{1}{32} \frac{\frac{T_1-x}{2}}{f(T_1)}.$$ Consequently I can infer that 
\begin{equation}\label{eq:negativecase2}
    \begin{aligned}
    \log \left( \frac{\partial S_2}{\partial y} \right) &\geq \log \left( 1+ \frac{\partial}{\partial y}\tilde{y}(x,T_1,y) + \left( \frac{\partial T_2}{\partial y}-1 \right) \left(\frac12 - \frac{\tilde{y}(x,T_1,y)-\frac{T_1-x}{4}}{\frac{T_1-x}{2}} \right) \right) \\
    & = \log \left( 1+ \frac{\partial}{\partial y}\tilde{y}(x,T_1,y)\right) + \log \left( 1  + \left( \frac{\partial T_2}{\partial y}-1 \right) \left(\frac12 + \tilde \delta \right)\right)  \\
    & \geq \log \left( 1+ \frac{\partial}{\partial y}\tilde{y}(x,T_1,y)\right)+ \log(1) + \log \left( \frac{\partial T_2}{\partial y} \right)-C \frac{\left(\frac{T_1-x}{2}\right)^2}{f(T_1)^2}
\end{aligned}
\end{equation}
where the last passage follows from Lemma \ref{lem:est2}. Finally it is possible to prove \eqref{eq:condition}, putting together \eqref{eq:negativecase1} with \eqref{eq:negativecase2}, applying Corollary \ref{cor:Kconve} and taking $K$ sufficiently large.
\end{pr}

As I did in section \ref{sec:jacobi}, I exploit the local nature of Jacobi equation to improve the last result. The following result is an easy consequence of Theorem \ref{thm:CDTheExample} and Corollary \ref{cor:jacobi}, and it will be useful in the end of this work.

\begin{corollary}\label{cor:entconvwithf}
Given two absolutely continuous measures $\mu_0,\mu_1\in \Prob(X)$, assume that there exists a map $T$ such that $T_\#\mu_0=\mu_1$, satisfying all the properties of Proposition \ref{prop:map}. Then, calling $M$ the midpoint selection presented in section \ref{sec:defimidpoint}, it holds 
\begin{align*}
    \Ent\big([M\circ (\id,T)]_\# (f\mu_0)\big) \leq \frac12 \Ent(f\mu_0) + \frac12 \Ent(T_\#(f\mu_0)),
\end{align*}
for every bounded measurable function $f:X\to \R^+$ with $\int f \de \mu_0=1$.
\end{corollary}

\section{Conclusions}

In this last section I conclude all the work done up to now. In particular I am going to show why this example is relevant, asking also to some open question related to strict curvature dimension bounds. First of all let me prove the most important result, which had already been anticipated in previous sections.  

\begin{theorem}\label{thm:CDcondatthelimit}
For suitable $k$ and $K$ the metric measure space $(X_{k,0},\di_\infty, \mathfrak{m}_{k,K,0})$ is a $\CD(0,\infty)$ space.
\end{theorem}

\begin{pr}
I am going to prove that, for every sequence of positive real numbers $(\varepsilon_n)_{n\in \setN}$ converging to zero, the sequence of metric measure spaces $(X_{k,\varepsilon_n},\di_\infty,\m_{k,K,\varepsilon_n})$ measured Gromov Hausdorff converges to $(X_{k,0},\di_\infty, \mathfrak{m}_{k,K,0})$. According to Theorem \ref{thm:stabilitycompact} and Theorem \ref{thm:CDTheExample} this is sufficient to conclude the proof, up to choose suitable $k$ and $K$. 

Define the function $f_n: X_{k,\varepsilon_n}\to X_{k,0}$ as 
\begin{equation*}
f_n(x,y)= \bigg(x,y \cdot \frac{f_{k,0}(x)}{f_{k,\varepsilon_n}(x)}\bigg) ,
\end{equation*}
it is immediate to notice that its image is actually $X_{k,0}$. Moreover its is easy to prove that 
\begin{equation*}
    (f_n)_\# \m_{k,K,\varepsilon_n} = \m_{k,K,0} \quad \text{and}\quad \big | \di_\infty \big( f_n(x_1,y_1), f_n(x_2,y_2) \big) - \di_\infty \big((x_1,y_1),(x_2,y_2)\big)\big| \leq 2 \varepsilon_n,
\end{equation*}
and this shows the desired measured Gromov-Hausdorff convergence.
\end{pr}

First of all notice that the space $(X_{k,0},\di_\infty, \mathfrak{m}_{k,K,0})$ has different topological dimensions at different regions of the space. In particular, this shows the non-constancy of topological dimension also for CD spaces, extending one of the results by Ketterer and Rajala \cite{ketterer2014failure}.
Furthermore, the space $(X_{k,0},\di_\infty, \mathfrak{m}_{k,K,0})$ is not a very strict $\CD(K,\infty)$ space for every $K\in\R$, in fact it is not weakly essentially non-branching (see Theorem \ref{thm:schultz}). In order to see this, it is sufficient to consider an absolutely continuous measure $\mu_0$ supported on $L$ and an absolutely continuous measure $\mu_1$ supported on $C$, and subsequently notice that every $\eta\in\OptGeo(\mu_0,\mu_1)$ is supported in branching geodesics. It is then possible to conclude that the weak CD condition is not sufficient to ensure any type of non-branching condition. Observe also that every $\eta\in\OptGeo(\mu_0,\mu_1)$ is not induced by a map, consequently the (weak) CD condition is not sufficient to ensure the existence of an optimal transport map, between two absolutely continuous marginals. Finally notice that the space $(X_{k,0},\di_\infty, \mathfrak{m}_{k,K,0})$ is an example of (weak) CD space which is not a very strict CD space, and this shows that this two notions of curvature dimension bounds are actually different.

For the last part of this section I need to introduce another type of curvature bounds, called strict CD condition, which is stronger than the weak CD condition, but is weaker than the very strict one.

\begin{definition}
A metric measure space $(X,\mathsf{d},\mathfrak{m})$ is called a strict $\CD(K,\infty)$ space
if for every absolutely continuous measures $\mu_0,\mu_1\in\ProbTwo(X)$
there exists an optimal geodesic plan $\eta\in \OptGeo(\mu_0,\mu_1)$, so that the entropy functional $\Ent$ satisfies the K-convexity inequality along $f\eta$ for every bounded measurable function $f : \Geo(X) \to \setR^+$ with $\int f \de \eta=1$.
\end{definition}

I am now going to prove that, for suitable constants, the spaces $(X_{k,\varepsilon},\di_\infty, \mathfrak{m}_{k,K,\varepsilon})$ with $0<\varepsilon<k$ are strict CD spaces, while their measured Gromov Hausdorff limit as $\varepsilon\to 0$ $(X_{k,0},\di_\infty, \mathfrak{m}_{k,K,0})$ is not. As a consequence the following proposition holds.

\begin{prop}\label{prop:strictCD}
The strict CD condition is not stable under measured Gromov Hausdorff convergence. 
\end{prop}

\noindent Before going into the details of the proofs, I want to make some clarifications. The fact that the spaces $(X_{k,\varepsilon},\di_\infty, \mathfrak{m}_{k,K,\varepsilon})$ are strict CD is a consequence of Corollary \ref{cor:entconvwithf} and of an iteration argument. On the other hand the space $(X_{k,0},\di_\infty, \mathfrak{m}_{k,K,0})$ cannot be a strict CD space, because of its particular topological structure I have already highlighted. 

\begin{prop}\label{prop:strictCDtrue}
For suitable $k$ and $K$ and every $0<\varepsilon<k$ the metric measure space $(X_{k,\varepsilon},\di_\infty, \mathfrak{m}_{k,K,\varepsilon})$ is a strict $\CD(0,\infty)$ space.
\end{prop}

\begin{pr}
For every $n\in \setN$ I am going to define a measurable map $G_n:X \to \Geo (X)$ by induction. In particular I introduce $G_0:X \to \Geo (X)$ as any measurable map such that $(e_0,e_1)\circ G_0(x)=(\id,T)$ $\mu_0$-almost everywhere, consequently $(G_0)_\# \mu_0 \in \OptGeo(\mu_0,\mu_1)$. Given $G_n:X \to \Geo (X)$, define $G_{n+1}:X \to \Geo (X)$ by imposing that:
\begin{enumerate}
    \item $e_r \circ G_{n+1}= e_r \circ G_{n}$ for every $r \in \big\{ \frac{k}{2^n}, k=0,\dots,2^n\big\}$
    \item $e_\frac{2k+1}{2^{n+1}} \circ G_{n+1}= M \big(e_\frac{k}{2^{n}}\circ G_n ,e_\frac{k+1}{2^{n}}\circ G_n\big)$ $\mu_0$-almost everywhere, where $M$ is the midpoint map defined as in Section \ref{sec:defimidpoint}.
\end{enumerate}
Notice that, if the optimal transport map that induces $\big(e_\frac{k}{2^{n}}\circ G_n ,e_\frac{k+1}{2^{n}}\circ G_n\big)_\# \mu_0$ satisfies all the properties of Proposition \ref{prop:map}, then also the maps that induce $\big(e_\frac{2k}{2^{n+1}}\circ G_{n+1} ,e_\frac{2k+1}{2^{n+1}}\circ G_{n+1}\big)_\# \mu_0$ and $\big(e_\frac{2k+1}{2^{n+1}}\circ G_{n+1} ,e_\frac{2k+2}{2^{n+1}}\circ G_{n+1}\big)_\# \mu_0$ satisfy all the properties of Proposition \ref{prop:map}. The reader can easily realize that this is a quite straightforward consequence of the definition of map $M$ and of its properties highlighted in Section \ref{section:CDproof}. This observation shows that the inductive procedure I have introduced can be done in accordance with the previous section, moreover it is possible to apply Corollary \ref{cor:entconvwithf} and obtain that
\begin{equation}\label{eq:convexityinduction}
\begin{split}
     \Ent\big((e_s)_\# [(G_{n+1})_\# (f\mu_0)] \big)&\leq \frac12\Ent\big((e_r)_\# [(G_{n+1})_\# (f\mu_0)]\big)+ \frac12\Ent\big((e_t)_\# [(G_{n+1})_\# (f\mu_0)]\big) \\
     &\quad - \frac K8 W_2^2\big((e_r)_\# [(G_{n+1})_\# (f\mu_0)],(e_t)_\# [(G_{n+1})_\# (f\mu_0)]\big),
     \end{split}
\end{equation}
where $r=\frac{2k}{2^{n+1}}$, $s=\frac{2k+1}{2^{n+1}}$, $t=\frac{2k+2}{2^{n+1}}$ and $f$ is any bounded measurable function with $\int f \de \mu_0=1$. Notice that, in order to infer \eqref{eq:convexityinduction}, I also used that the map $e_r \circ G_n$ is injective outside a $\mu_0$-null set, as a consequence of Proposition \ref{prop:aeinjectivity}.\\
An inductive argument allows to conclude that for every $n\in \setN$ it holds
\begin{equation*}
    \Ent\big((e_r)_\# [(G_n)_\# (f\mu_0)] \big)\leq (1-r)\Ent(f\mu_0) + r\Ent(T_\#(f\mu_0)) - \frac K2 r(1-r) W_2^2(f\mu_0,T_\#(f\mu_0)),
\end{equation*}
for every $r \in \big\{ \frac{k}{2^n}, k=0,\dots,2^n\big\}$, and every bounded measurable function $f$ with $\int f \de \mu_0=1$. In fact, this is completely obvious for $n=0$, and assuming it true for an $n$ it is possible to easily deduce it for $n+1$, using \eqref{eq:convexityinduction}.\\
It is now easy to notice that the first property in the definition of $G_{n+1}$ given $G_n$, ensures the existence of a measurable map $G:X\to \Geo(X)$, such that $G_n\to G$ uniformly. Furthermore it is obvious that 
\begin{equation*}
    (e_r)_\# [(G_n)_\# (f\mu_0)]=(e_r)_\# [G_\# (f\mu_0)]
\end{equation*}
for every $r \in \big\{ \frac{k}{2^n}, k=0,\dots,2^n\big\}$, and every bounded measurable function $f$ with $\int f \de \mu_0=1$. Consequently it holds that
\begin{equation*}
    \Ent\big((e_r)_\# [G_\# (f\mu_0)] \big)\leq (1-r)\Ent((f\mu_0)) + r\Ent(T_\#(f\mu_0)) - \frac K2 r(1-r) W_2^2(f\mu_0,T_\#(f\mu_0)),
\end{equation*}
for every dyadic time $r$ and every suitable function $f$.
Then the lower semicontinuity of $\Ent$ allows to infer that the $K$-convexity inequality of the entropy is satisfied along $G_\#(f \mu_0)$ for every suitable function $f$.
Finally I can conclude by observing that every optimal geodesic plan of the type $F \cdot G_\# \mu_0$, for a measurable function $F:\Geo(X)\to \R^+$ with $\int F \de (G_\# \mu_0)=1$, can be written as $G_\#(f \mu_0)$ for a suitable measurable $f$ with $\int f \de \mu_0=1$, since $G$ is clearly injective.
\end{pr}

\begin{prop}\label{prop:strictCDnottrue}
For every $k$ and $K$ the metric measure space $(X_{k,0},\di_\infty, \mathfrak{m}_{k,K,0})$ is not a strict $\CD(0,\infty)$ space.
\end{prop}

\begin{figure}
\begin{center}

\begin{tikzpicture} 

\shade[left color=black!22!white,right color=white] (0,0)--(6,0)--(6,3)--(1,0.5)..controls (0,0)..(-1,0)--cycle;
\draw[thick](-6,0)--(6,0);
%\filldraw[black] (0,0) circle (1pt);
%\draw[thick,dotted](0,0)--(7,3.5);
\draw[thick](6,0)-- (6,3)--(1,0.5);
\draw[thick](-1,0)..controls (0,0)..(1,0.5); 
\filldraw[white] (-7,0) circle (1pt);
%\node at (3,1.6)[label=north:$f_{k,0}(x)$] {};
\draw[ultra thick](-3,0)--(-1.5,0);
\filldraw[black] (-1.5,0) circle (1pt);
\filldraw[black] (-3,0) circle (1pt);
\node at (-2.25,0)[label=south:$\mu_0$] {};
\shade [left color=black!37!white,right color=black!22!white] (4.5,0)--(3,0)-- (3,1.5)--(4.5,2.25)--cycle;
\draw [ultra thick] (4.5,0)--(3,0)-- (3,1.5)--(4.5,2.25)--cycle;
\node at (3.75,0)[label=south:$\mu_1$] {};
\shade [left color=black!47!white,right color=black!40!white](-0.75,0)--(0.75,0)-- (0.75,0.375)..  controls (-0.1,0)..(-0.25,0)--cycle;
\draw [ultra thick] (-0.75,0)--(0.75,0)-- (0.75,0.375)..  controls (-0.1,0)..(-0.25,0)--cycle;
\node at (0,0)[label=south:$\mu_s$] {};

\end{tikzpicture}

\caption{A representation of the geodesic $(\mu_t)_{t\in[0,1]}$.}
\label{fig:notstrict}
\end{center}
\end{figure}
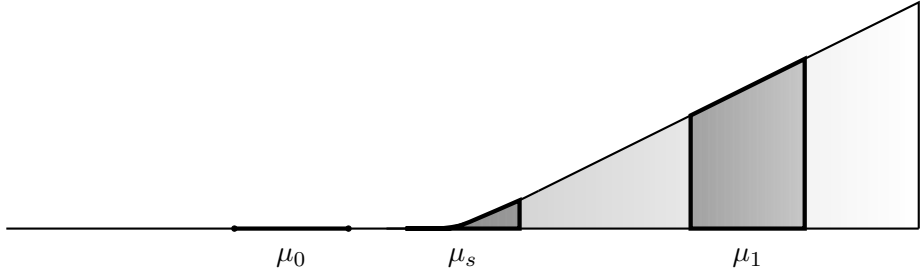

\begin{pr}
In this proof I denote with $\m$ the measure $\mathfrak{m}_{k,K,0}$, in order to simplify the notation.  \\
For every $t\in [0,1]$ define the measure 
\begin{equation*}
    \mu_t = \frac{1}{\mathfrak{m}\big([-\frac12 +t,-\frac14+t]\times \R\big)}\cdot \restr{\mathfrak{m}}{[-\frac12 +t,-\frac14+t]\times \R} =  \frac{4}{C_{K}} \cdot \restr{\mathfrak{m}}{[-\frac12 +t,-\frac14+t]\times \R},
\end{equation*}
see Figure \ref{fig:notstrict} in order to visualize it. It is easy to realize that $(\mu_t)_{t\in[0,1]}$ is the unique geodesic connecting $\mu_0$ and $\mu_1$, along which the entropy functional is convex. Moreover let $\eta\in \Prob\big(C([0,1]),X\big)$ such that $(e_t)_\# \eta = \mu_t$ for every $t\in[0,1]$, I am going to prove that the entropy functional is not convex along $F\eta$, for a suitable bounded measurable function $F:\Geo(X)\to \R^+$ with $\int F \de \eta=1$. Before going on, let me point out that every Wasserstein geodesic in $\Prob\big(C([0,1]),X\big)$ which connects a measure on $L$ to a measure on $C$ (and thus $\eta$ in particular), consists only of ``horizontal" transports. Therefore the only useful coordinate, in order to evaluate the distance $\di_\infty$, will be the $x$ coordinate. As a consequence every such optimal geodesic plan (and $\eta$ in particular) will only depend on the $x$ coordinate. Some of the considerations I will do in this proof actually follows from this observation.\\
Now, define the set 
\begin{equation*}
    A:= \{(x,y)\in \R^2 : f_{k,0}(x)>0 \text{ and } 0\leq 2y \leq f_{k,0}(x)  \}\subset X_{k,0},
\end{equation*}
and the quantity
\begin{equation*}
    C'_{K}= \int_0^\frac12 e^{-K y^2} \de y.
\end{equation*}
Then fix a time $\bar t$ such that $\mu_{\bar t}$ is concentrated in $C$, consider the map $\tilde F:\Geo(X)\to \R^+$ defined as $\tilde F:= \frac{C_K}{C'_K} \cdot \chi_A \circ e_{\bar t}$ and call $\tilde \mu_1 = (e_1)_\# (\tilde F \eta)$. Notice that $\int \tilde F \de \eta=1$, thus $\tilde \mu_1$ is a probability measure and it is absolutely continuous with respect to $\m$, with density $\tilde \rho_1$ bounded above by $\frac{4}{C'_K}$, as a consequence $\m(\{\tilde \rho_1 >0\})\geq \frac{C'_K}{4}$. Now, suppose that $\m(\{\tilde \rho_1 >0\})= \frac{C'_K}{4}$, then $\tilde \rho_1 \equiv \frac{4}{C'_K}$ on $\{\tilde \rho_1 >0\}$ and therefore 
\begin{equation*}
    \Ent(\tilde \mu_1)= \log \left( \frac{4}{C'_K} \right) = \Ent \big((e_{\bar t})_\# (\tilde F \eta)\big).
\end{equation*}
On the other hand 
\begin{equation*}
    \Ent \big((e_{0})_\# (\tilde F \eta)\big) = \Ent(\mu_0)  = \log \left( \frac{4}{C_K} \right)<\log \left( \frac{4}{C'_K} \right),
\end{equation*}
and consequently the entropy functional is not convex along $\tilde F \eta$.\\
Otherwise, suppose that $\m(\{\tilde \rho_1 >0\})> \frac{C'_K}{4}$, call $S:=\{\tilde \rho_1 >0\}$ and define the set $S_x:=\{(x',y')\in S : x'=x\}$, for every $x\in [-1,1]$. Now consider $\m_1:=(\p_1)_\# \m$ and denote by $(\m_x)_{x\in[-1,1]}\subset \Prob(\R)$ the disintegration of $\m$ with respect to the projection map $\p_1$. Notice that, since $\eta$ depends only on the $x$ coordinate, then $(\p_1)_\#\tilde\mu_1= \frac{4}{C_K} \cdot \restr{\m_1}{\left[\frac 12,\frac 34\right]}$. Moreover, since the density $\tilde \rho_1$ is bounded above by $\frac{4}{C'_K}$, it holds that $\m_x(S_x)\geq \frac{C'_K}{C_K}$, for $\m_1$-almost every $x\in \left[\frac 12,\frac 34\right]$. Furthermore, the assumption on $S$, that is $\m(S)> \frac{C'_K}{4}$, ensures that $\m_x(S_x)> \frac{C'_K}{C_K}$ for a $\restr{\m_1}{\left[\frac 12,\frac 34\right]}$-positive set of $x$, therefore
\begin{equation}\label{eq:integralest}
    \int_\frac 12 ^\frac34 \log \left(  \m_x(S_x)\right)\de\m_1(x)  > \frac{C_K}{4}\log \left(\frac{C'_K}{C_K}\right) .
\end{equation}
On the other hand, for every positive constant $c>0$ define the set $S^c:=\{\tilde\rho_1>c\}$ and call $S^c_x:=\{(x',y')\in S^c : x'=x\}$ for every $x\in[-1,1]$. Notice that, for every sufficiently small constant $c$, since $\tilde \rho_1$ is bounded and $(\p_1)_\#\tilde\mu_1= \frac{4}{C_K} \cdot \restr{\m_1}{\left[\frac 12,\frac 34\right]}$, the quantity $\m_x(S_x)$ is uniformly bounded from below for $\m_1$-almost every $x\in\left[\frac 12,\frac 34\right]$. Consequently, it is possible to apply the monotone convergence theorem and deduce that there exists a constant $\bar c>0$ such that 
\begin{equation*}
    \int_\frac 12 ^\frac34 \log \left(  \m_x(S^{\bar c}_x)\right)\de\m_1(x)  > \frac{C_K}{4}\log \left(\frac{C'_K}{C_K}\right) .
\end{equation*}
Now, define the measurable map $ F:\Geo(X)\to \R^+$
\begin{equation*}
    F :=  \frac{4}{C_K} \tilde F \cdot \left(\frac{1}{\m_x(S^{\bar c}_x)}\frac{\chi_{S^{\bar c}}(x,y)}{\tilde\rho_1(x,y)} \circ e_1 \right).
\end{equation*} 
I have already noticed that the quantity $\m_x(S^{\bar c}_x)$ is uniformly bounded from below for $\m_1$-almost every $x\in\left[\frac 12,\frac 34\right]$, thus $F$ is well defined and bounded, moreover it holds that
\begin{align*}
    \int F \de \eta &= \int \frac{4}{C_K} \cdot \left(\frac{1}{\m_x(S^{\bar c}_x)}\frac{\chi_{S^{\bar c}}(x,y)}{\tilde\rho_1(x,y)} \circ e_1 \right) \de \tilde F \eta = \frac{4}{C_K}\int  \frac{1}{\m_x(S^{\bar c}_x)}\frac{\chi_{S^{\bar c}}(x,y)}{\tilde \rho_1(x,y)}   \de \tilde\mu_1 (x,y)\\
    &= \frac{4}{C_K}\int  \frac{\chi_{S^{\bar c}}(x,y)}{\m_x(S^{\bar c}_x)}   \de \m (x,y) 
     = \frac{4}{C_K} \int_\frac12 ^\frac 34 \int \frac{\chi_{S^{\bar c}}(x,y)}{\m_x(S^{\bar c}_x)} \de \m_x(y) \de \m_1(x)
    =  \frac{4}{C_K} \int_\frac12 ^\frac 34 \de \m_1(x)=1.
\end{align*}
In particular, observe that a computation similar to this last one shows that
\begin{equation*}
    (e_{1})_\# ( F \eta) = \frac{4}{C_K} \frac{\chi_{S^{\bar c}}(x,y)}{\m_x(S^{\bar c}_x)} \cdot \m,
\end{equation*}
thus it is possible to estimate its entropy:
\begin{align*}
    \Ent\big((e_{1})_\# ( F \eta)\big) &= \int  \frac{4}{C_K} \frac{\chi_{S^{\bar c}}(x,y)}{\m_x(S^{\bar c}_x)}\log \left(  \frac{4}{C_K} \frac{\chi_{S^{\bar c}}(x,y)}{\m_x(S^{\bar c}_x)}\right) \de \m  \\
    &=  \int_\frac12 ^\frac 34 \int \frac{4}{C_K} \frac{\chi_{S^{\bar c}}(x,y)}{\m_x(S^{\bar c}_x)}\log \left(  \frac{4}{C_K} \frac{\chi_{S^{\bar c}}(x,y)}{\m_x(S^{\bar c}_x)}\right) \de \m_x(y) \de \m_1(x)\\
      &=  \int_\frac12 ^\frac 34 \int_{S^{\bar c}_x} \frac{4}{C_K} \frac{1}{\m_x(S^{\bar c}_x)}\log \left(  \frac{4}{C_K} \frac{1}{\m_x(S^{\bar c}_x)}\right) \de \m_x(y) \de \m_1(x)\\
      &= \frac{4}{C_K} \int_\frac 12 ^\frac34 \log \left(  \frac{4}{C_K} \frac{1}{\m_x(S^{\bar c}_x)}\right)\de \m_1(x) \\
      &= \log \left(  \frac{4}{C_K} \right)+ \frac{4}{C_K} \int_\frac 12 ^\frac34 -\log \left(  \m_x(S^{\bar c}_x)\right)\de \m_1(x) < \log \left( \frac{4}{C'_K} \right),
\end{align*}
where the last inequality follows from \eqref{eq:integralest}.
On the other hand $(e_{\bar t})_\#(F\eta)\ll (e_{\bar t})_\# (\tilde F \eta)$ and consequently Jensen's inequality ensures that
\begin{equation*}
    \Ent\big((e_{\bar t})_\#(F\eta)\big) \geq  \log \left( \frac{4}{C'_K} \right),
\end{equation*}
Furthermore, it is easy to realize that 
\begin{equation*}
    (\p_1)_\#\big[(e_{1})_\# ( F \eta)\big]= \frac{4}{C_K} \cdot \restr{\m_1}{\left[\frac 12,\frac 34\right]}
\end{equation*}
and thus, as before, I have
\begin{equation*}
    \Ent \big((e_{0})_\# ( F \eta)\big) = \Ent(\mu_0)  = \log \left( \frac{4}{C_K} \right)<\log \left( \frac{4}{C'_K} \right).
\end{equation*}
Putting together this last three inequalities it is easy to realize that the entropy functional is not convex along $F \eta$.
\end{pr}

\noindent Notice that this last result shows in particular that the strict CD condition and the weak one are two actually different notions. Moreover, the combination of Proposition \ref{prop:strictCDtrue} and \ref{prop:strictCDnottrue} obviously yields Proposition \ref{prop:strictCD}, according to the proof of Theorem \ref{thm:CDcondatthelimit}. On the other hand, observe that this work does not allow to disprove the stability of the very strict CD condition. In fact the proof of Proposition \ref{prop:strictCDtrue} heavily relies on an approximation argument, thus it is impossible to modify it, in order to prove the very strict CD condition for the spaces $(X_{k,\varepsilon},\di_\infty, \mathfrak{m}_{k,K,\varepsilon})$. However, in my opinion, this example leaves little hope for the very strict CD condition to be stable.\\

\noindent {\scshape Aknowlegments} : This article contains part of the work I did for my master thesis, that was supervised by Luigi Ambrosio and Karl-Theodor Sturm.

\bibliography{example}

\begin{thebibliography}{10}

\bibitem{ambrosio2005gradient}
L.~Ambrosio, N.~Gigli, and G.~Savare.
\newblock {\em Gradient Flows: In Metric Spaces and in the Space of Probability
  Measures}.
\newblock Lectures in Mathematics. ETH Z{\"u}rich. Birkh{\"a}user Basel, 2005.

\bibitem{Gigli_2015}
N.~Gigli, A.~Mondino, and G.~Savaré.
\newblock Convergence of pointed non-compact metric measure spaces and
  stability of {R}icci curvature bounds and heat flows.
\newblock {\em Proceedings of the London Mathematical Society}, 111:1071--1129,
  2015.

\bibitem{ketterer2014failure}
C.~Ketterer and T.~Rajala.
\newblock Failure of topological rigidity results for the measure contraction
  property.
\newblock {\em Potential Analysis}, 42, 2014.

\bibitem{lottvillani}
J.~Lott and C.~Villani.
\newblock Ricci curvature for metric-measure spaces via optimal transport.
\newblock {\em Annals of Mathematics}, 169:903--991, 2009.

\bibitem{rajala2013failure}
T.~Rajala.
\newblock Failure of the local-to-global property for {$\CD (K, N)$} spaces.
\newblock {\em Ann. Sc. Norm. Super. Pisa Cl. Sci.}, 15:45--68, 2016.

\bibitem{rajalasturm}
T.~Rajala and K.-T. Sturm.
\newblock Non-branching geodesics and optimal maps in strong {$\CD(K,\infty)$}
  spaces.
\newblock {\em Calculus of Variations and Partial Differential Equations},
  50:831--846, 2014.

\bibitem{schultz2017existence}
T.~Schultz.
\newblock Existence of optimal transport maps in very strict {$\CD(K,\infty)$}
  spaces.
\newblock {\em Calculus of Variations and Partial Differential Equations}, 57,
  2018.

\bibitem{Schultz2019EquivalentDO}
T.~Schultz.
\newblock Equivalent definitions of very strict {$\CD(K,N)$} spaces.
\newblock {\em arXiv preprint}, 2019.

\bibitem{Schultz2019OnOO}
T.~Schultz.
\newblock On one-dimensionality of metric measure spaces.
\newblock {\em Proc. Amer. Math. Soc.}, 149:383--396, 2020.

\bibitem{sturm2006}
K.-T. Sturm.
\newblock On the geometry of metric measure spaces.
\newblock {\em Acta Math.}, 196(1):65--131, 2006.

\bibitem{sturm2006ii}
K.-T. Sturm.
\newblock On the geometry of metric measure spaces. {II}.
\newblock {\em Acta Math.}, 196(1):133--177, 2006.

\bibitem{villani2008}
C.~Villani.
\newblock {\em Optimal transport -- Old and new}.
\newblock Grundlehren der mathematischen Wissenschaften. Springer, 2008.

\end{thebibliography}

\bibliographystyle{abbrv} 

\end{document}